\magnification1150
\input epsf.sty
\input  amssym.tex
\input color

\def\sqr#1#2{{\vcenter{\hrule height.#2pt              
     \hbox{\vrule width.#2pt height#1pt\kern#1pt
     \vrule width.#2pt}
     \hrule height.#2pt}}}
\def\square{\mathchoice\sqr{5.5}4\sqr{5.0}4\sqr{4.8}3\sqr{4.8}3}
\def\qed{\hskip4pt plus1fill\ $\square$\par\medbreak}


\centerline{\bf Dynamics of (Pseudo) Automorphisms of 3-space: }

\centerline{\bf  Periodicity versus positive entropy}

\bigskip
\centerline{Eric Bedford\footnote*{Supported in part by the NSF}  and Kyounghee Kim\footnote{\dag} {Supported in part by the NSA} }

\bigskip 

\bigskip

\noindent{\bf \S0.  Introduction.}   We consider the family of birational maps of 3-space which may be written in affine coordinates as
$$f_{\alpha,\beta}: (x_1,x_2,x_3)\mapsto \left(x_2,x_3,{\alpha_0 +\alpha_1 x_1+ \alpha_2 x_2 + \alpha_3x_3\over \beta_0+\beta_1 x_1+\beta_2 x_2 + \beta_3 x_3}\right).\eqno(0.1) $$
The algebraic iterates $f_{\alpha,\beta}^n:=f_{\alpha,\beta}\circ\cdots\circ f_{\alpha,\beta}$ are rational maps for all $n\in{\bf Z}$.  Here we study the dynamics of $f=f_{\alpha,\beta}$, by which we mean the behavior of $f^n$ as $n\to\pm\infty$.    We have invertible dynamics since $f$ has a rational inverse, but  it does not behave like a diffeomorphism (or even a homeomorphism).  There are two difficulties if we want to regard $f$ as a mapping of points.  First, there is the set of indeterminacy ${\cal I}(f)$; $f$ blows up each point of ${\cal I}(f)$ to a variety of positive dimension.   Second, there can be hypersurfaces $E$ which are exceptional, in the sense that the codimension of $f(E-{\cal I}(f))$ is at least 2.  We will say that $f$ is a pseudo-automorphism if neither $f$ nor $f^{-1}$ has an exceptional hypersurface.  In dimension 2, every pseudo-automorphism is in fact an automorphism.   However,  for pseudo-automorphisms,  indeterminate behaviors are possible in higher dimension which have no analogue in dimension~2.

Given a rational map $f:X\dasharrow X$ there is a well-defined pullback map  on cohomology,  $f^*:H^*(X)\to H^*(X)$.      Passage to cohomology, however,  may not be compatible with iteration because the identity  $(f^*)^n=(f^n)^*$ may not be valid.   Given a birational map $f$ in dimension 2, Diller and Favre [DiF] showed that there is a new manifold $\pi:Y\to X$ such that the iterates of the induced map $f_Y$ behave naturally on cohomology, in the sense that  $(f_Y^*)^n=(f_Y^n)^*$.  In dimension greater than 2, however, no such theorem is known.

Given a rational map of ${\bf P}^n$ we may consider modifications $\pi:X\to{\bf P}^n$, where $\pi$ is a morphism which is birational.  This induces a rational map $f_X:=\pi^{-1}\circ f\circ \pi$ of $X$, which might have pointwise properties which are different from those of the original $f$. If $f_X$ is a pseudo-automorphism, then $f_X$ acts naturally on $H^{1,1}(X)$.   The exponential rate of growth of $f^n$ on $H^{p,p}$:  $\delta_p(f): = \lim_{n\to\infty}||f^{n*}|_{H^{p,p}(X)}|| ^{1/n}$ is known as the $p$th dynamical degree and is a birational invariant (see [DS]).

Within the family (0.1) we find the first known examples of pseudo-automorphisms of positive entropy on blowups of ${\bf P}^3$:
\proclaim Theorem 1.  Suppose that $\alpha=(a, 0, \omega, 1)$ and $\beta=(0,1,0,0)$ where $a \in {\bf C} \setminus\{0\}$ and $\omega$ is a non-real cube root of the unity.   Then  there is a modification $\pi:Z\to{\bf P}^3$ such that $f_Z$ is a pseudo-automorphism.
The dynamical degrees $\delta_1(f)=\delta_2(f)\cong 1.28064>1$ are equal and are given by the largest root of $t^8-t^5-t^4-t^3+1$.  The entropy of $f_Z$ is the logarithm of the dynamical degree and is thus positive.

\proclaim Theorem 2. For the mappings in Theorem 1, there is a 1-parameter family of surfaces $S_c\subset Z$, $c\in {\bf C}$ which have the invariance $fS_c=S_{\omega c}$.   For generic $c$, $S_c$ is $K3$, and the restriction $f^3|_{S_c}$ is an automorphism.  For generic $c$ and $c'$, the surfaces $S_c$ and $S_{c'}$ are biholomorphically inequivalent, and the automorphisms $f^3|_{S_c}$ and $f^3|_{S_{c'}}$ are not smoothly conjugate.   

The surface $S_0$ is invariant, and the restriction  $f_{S_0}$ is an automorphism which has the same entropy as $f$.  This is smaller than the entropy of the automorphism constructed in [M2, Theorem 1.2] and is thus the smallest known entropy for a {\it projective} $K3$ surface automorphism.

Closely related to the dynamics of $f_Z$ is the (1,1)-current $T^+$ which is expanded by $f_Z^*$, and a current $T^-$ for $f^{-1}_Z$.  This is obtained in \S7, as well as the invariant (2,2)-current $T^+\wedge T^-$.  The slices of $T^\pm$ and $T^+\wedge T^-$ on the surfaces $S_c$ give the expanded/contracted currents, as well as the unique invariant measure, for the automorphism $f|_{S_c}$.

The following mappings have quadratic degree growth and complete integrability:
 
\proclaim Theorem 3. Suppose that $\beta=(0,1,0,0)$ and either $\alpha =(0,0,\omega,1)$ or $\alpha=(a,0,1,1)$ where $a \in {\bf C} \setminus\{1\}$, $\omega\ne1$, and $\omega^3=1$.  Then the degree of $f^n$ grows quadratically in $n$.  Further, there is a modification $\pi:Z\to{\bf P}^3$ such that $f_Z$ is a pseudo-automorphism.  There is a two-parameter family of surfaces $S_c$, $c=(c_1,c_2)\in{\bf C}^2$ which are invariant under $f^3$.  For generic $c$ and $c'$, $S_c$ is a smooth $K3$ surface, and $S_{c}\cap S_{c'}$ is a smooth elliptic curve.

For the mappings in Theorems 1 and 3, $f$ is reversible on the level of cohomology:  $f_Z^*$ is conjugate to $(f_Z^{-1})^*=(f_Z^*)^{-1}$.  The identity $\delta_1(f)=\delta_2(f)$ for such maps is a consequence of the duality between $H^{1,1}$ and $H^{2,2}$, so they are not cohomologically hyperbolic, in the terminology of [G2].  For each of these maps, the family of invariant $K3$ surfaces becomes singular at an invariant 8-cycle ${\cal R}$ of rational surfaces (see (7.2)).  We show that the restriction $f|_{\cal R}$ is not birationally conjugate to a surface automorphism: see Appendix C for the maps in Theorem 1 and Proposition 8.2 for the maps in Theorem 3.  By Corollary 1.6, then, we have:  
\proclaim Theorem 4.  Let $f$ be a map from Theorems 1 and 3.  If $a\ne1$, then $f$ is not birationally conjugate to an automorphism.


We note that for birational surface maps, the degree growth of the iterates determines whether the map is birationally conjugate to an automorphism:  This occurs if and only if either  (i) the degrees are bounded or degree growth is quadratic (see [DiF]), or (ii) if the dynamical degree is a Salem number (see  [BC]).  Theorem 4 shows that this result does not hold in dimension 3.

We will also determine which mappings $f_{\alpha,\beta}$ are periodic, or finite order, in the sense that $f^p=id$ for some $p>0$.  In contrast to Theorem 4, it was shown by de Fernex and Ein [dFE] that if $f$ is a rational map of finite order, then there is a modification $f_X$ as above, which is an automorphism of $X$.  If $f_X$ is periodic, then $f_X^*$ will also be periodic.  

In (4.1) and (4.2) we identify conditions which are necessary for $f$ to be periodic and are sufficient for the existence of a space $Z=Z_{\alpha,\beta}$ such that $f_Z$ is a pseudo-automorphism.  We show that for a map in (0.1), if $f_Z^*$ is periodic, then $f$ also turns out to be periodic.  The birational map (0.1) may also be considered as a  3-step linear fractional recurrence: given $z_0,z_1,z_2$, we define a sequence $\{z_n\}$ by
$$z_{n+3} = { \alpha_0 + \alpha_1 z_{n} +\alpha_2 z_{n+1} + \alpha_3 z_{n+2} \over \beta_0 + \beta_1 z_n + \beta_2 z_{n+1} +\beta_3 z_{n+2}} .\eqno(0.2)$$
The recurrence (0.2) is said to be periodic if the sequence $\{z_n\}$ is periodic for all choices of initial terms $z_0$, $z_1$ and $z_2$. Equivalently, $f_{\alpha,\beta}^p=id$ for some $p$. For all $r>0$ there are $r$-step recurrences of the form (0.2).  In [BK2] we determined the possible periods for 2-step linear fractional recurrences.  McMullen [M1] has explained the periods that arise by showing that the corresponding (2-dimensional) $f_{\alpha,\beta}$ represent certain Coxeter elements.  

Here we determine all possible periods for 3-step recurrences (0.2).  To rule out trivial cases, we assume that the coefficients satisfy (2.3), and we have:
\proclaim Theorem 5.  The only nontrivial periods for (0.2) are 8 and 12.  Each periodic recurrence is equivalent to one of the following:
$$ z_{n+3}= {1+z_{n+1}+z_{n+2}\over z_n} \qquad z_{n+3} ={ -1 -z_{n+1} + z_{n+2}\over z_n} \eqno{\rm (period \ 8)}$$
$$  z_{n+3} =   {\eta/(1-\eta)+\eta z_{n+1}+z_{n+2} \over \eta^2 + z_n} \ \ \ \ \ \ \eta^3=-1 \eqno {\rm (period\ 12)}$$
In the notation of (0.1), the first case corresponds to $\beta=(0,1,0,0)$, $\alpha=(\pm 1, 0, \pm1,1)$, and the second case to $\beta=(\eta^2,1,0,0)$, $\alpha = (\eta/(1-\eta),0,\eta,1)$.

Each of these mappings has a different structure; these structures are described in Theorems 6.10 and 6.11. The first period 8 recurrence above was found by Lyness [L], and the  second one was found by Cs\"ornyei and Laczkovic [CL] (see also [CGM1]).  We note that the period 12 recurrences are the case $k=3$  of a general phenomenon exhibited in [BK4]:  {\sl For each $k$, there are $k$-step linear fractional recurrences with period $4k$.}  There is a literature dealing with $r$ step recurrences of the form (0.2). We refer to  the books [KoL], [KuL], [GL],  [CaL] and the extensive bibliographies they contain.   That direction of research  is largely concerned with the case where the structural parameters $\alpha$, $\beta$, as well as the dynamical points, are real and positive.  This avoids the difficulty that the denominator in $(0.2)$ might vanish, causing the expression to be undefined; but the restriction to positive numbers leads to a subdivision into a large number of distinct cases to be treated separately.  

In working with the family $f_{\alpha,\beta}$, we work with the pointwise iterates as much as possible, but this runs into difficulties if the orbit enters the indeterminacy locus.  We can often deal with this by blowing up certain subsets.  In this way we convert these subsets into hypersurfaces, and we then deal with the hypersurfaces by passing to $f^*$ on $Pic$.  This allows us to convert many difficulties with indeterminate orbits into more tractable problems of Linear Algebra.

This paper is organized as follows.  \S1 assembles some general information about rational maps and the geometry of blowing up.  \S2 gives the specific behaviors of the maps (0.1).  It is evident, then, that there are two possibilities, defined by (3.1), which we call ``critical'' and ``non-critical,'' and in \S3 we show that any  periodic map must be critical.  We study the structure of general critical maps in \S4.  In Theorem 5.1 we show that if $f$ is a critical map satisfying (5.1), then $f_{  Z}$ is a pseudo-automorphism.  Pseudo-automorphisms are discussed in \S5, together with the possibilities for the induced map $f^*_{  Z}$ on cohomology.  In \S6  we determine the periodic mappings and  give the proof of Theorem 5. In \S7 we give the proof of Theorems~1 and 2.  At the end of \S7 we present a different pseudo-automorphism with positive entropy; it has properties similar to those given in Theorems~1 and 2, but we do not discuss it in detail. The proof of Theorem 3 is given in \S8.

\bigskip
\noindent{\bf \S1.  Rational Maps }   A rational map $f:{\bf P}^d\dasharrow{\bf P}^d$ is given by a $(d+1)$-tuple of homogeneous polynomials, all of the same degree: $f=[f_0:\cdots:f_d]$.   We may divide $f$ by  ${\rm g.c.d.}(f_0,\dots,f_d)$ so that $f_i$'s have no common polynomial factor.   We define the degree of $f$, ${\rm deg}(f)$, to be the (common) degree of the $f_j$'s.  The {\it  indeterminacy locus } of $f$ is defined by
$${\cal I}(f)=\{x\in{\bf P}^d:f_0(x)=\cdots=f_d(x)=0\}$$
and is a subvariety of codimension at least 2, and  $f$ defines a holomorphic mapping $f:{\bf P}^d\setminus{\cal I}(f)\to {\bf P}^d$.  
 If $S$ is an irreducible subvariety of ${\bf P}^d$, and $S\not\subset{\cal I}(f)$, we define the {\it strict transform}, written simply as $f(S)$, to be the closure of $f(S-{\cal I}(f))$.
We say that an irreducible variety $V$ is {\it exceptional} for a rational mapping $f$ if $V\not\subset{\cal I}(f)$, and if the dimension of $f(V-{\cal I}(f))$ is strictly less than the dimension of $V$.  Following [DO, p.\ 64], we say that $f:X\dasharrow Y$ is a {\it pseudo-isomorphism} if $f$ is birational, and if neither $f$ nor $f^{-1}$ has an exceptional hypersurface.  It follows that if $f$ is a pseudo-isomorphism, then $f:X\setminus{\cal I}(f) \to Y\setminus{\cal I}(f^{-1})$ is biholomorphic.  If $X=Y$, we say that $f$ is a {\it pseudo-automorphism}.

\proclaim Theorem 1.1.  If $f:X\dasharrow Y$ is a pseudo-isomorphism between 3-dimensional manifolds, then the indeterminacy locus has no isolated points.

\noindent{\it Proof. }   Suppose that there is an isolated point $p\in{\cal I}(f)$.  Since $f^{-1}$ has no exceptional hypersurfaces, $f$ must blow $p$ up to a curve $C'\subset Y$.  Now we consider the behavior of $f^{-1}$ on $C'$.  We must have $C'\subset{\cal I}(f^{-1})$, for if $f^{-1}$ is regular at a point $q\in C'$, then $f^{-1}$ must map an open subset of $C'$ to $p$.  Thus the jacobian of $f^{-1}$ must vanish at $q$.  Since the jacobian vanishes on a hypersurface, $f^{-1}$ would have an exceptional hypersurface containing $q$.  Thus $q$ must be indeterminate.  Since the total transform of $q$ under $f^{-1}$ is  given by $\bigcap_{\epsilon>0}\overline{ \left(f^{-1}(B(q,\epsilon)-{\cal I}(f^{-1}))\right)}$, it must be connected, and it must be a curve $C$ containing $p$.  But since $p$ was an isolated point of ${\cal I}(f)$, there are nearby points $p'\in C-{\cal I}(f)$.  Since $f$ is regular at these points, it must map them to $q$, and thus $f$ must have an exceptional hypersurface.  By this contradiction, we see that ${\cal I}(f)$ has no isolated points.
\qed

For a rational map $f:{  X}\dasharrow{  X}$, we consider the iterates $f^j=f\circ\cdots\circ f$, $j>0$.  If $\Sigma$ is an irreducible hypersurface, then $\Sigma\not\subset{\cal I}(f^j)$ for reasons of dimension, so we may consider the sequence of varieties $V_j:=f^j(\Sigma)$, for $j>0$.  Since we will be interested in knowing to what extent the iterates of $f$ behave like a pointwise-defined dynamical system, we note:  {\sl   If $S\not\subset{\cal I}(g)$ is irreducible and if $g(S)\not\subset{\cal I}(f)$, then $S\not\subset {\cal I}(f\circ g)$, and $f(g(S))=(f\circ g)(S)$.}
We may also define $f$ at points of indeterminacy.  Let $\gamma_f=\{(x,y)\in({\bf P}^d-{\cal I})\times{\bf P}^d:y=f(x)\}$ denote the graph of $f$ at its regular points, and we let $\Gamma$ denote the closure of $\gamma_f$ inside ${\bf P}^d\times{\bf P}^d$.  It follows that $\Gamma$ is an irreducible variety of dimension $d$, and there are holomorphic projections $\pi_j:\Gamma\to{\bf P}^d$, $j=1,2$, onto the first and second factors, respectively, and we have $f=\pi_2\circ\pi_1^{-1}$ on ${\bf P}^d-{\cal I}$.    For a point $p\in {\bf P}^d$, we define the {\it total transform} to be $f_*p:=\pi_2(\pi_1^{-1}p)$, and then we define $f_*(S):=\bigcup_{p\in S}f_*p$.   It is easily seen that we have:  {\sl If $\Sigma$ is an irreducible hypersurface,  then $f_*(g(\Sigma))\supset(f\circ g)(\Sigma)$.}

\proclaim Proposition 1.2.  Suppose that $f:{  X}\dasharrow{  X}$ is rational, and suppose that for each exceptional hypersurface $E$ and for $m>0$, we have $f^m (E- {\cal I}) \not\subset {\cal I}$ .  If follows that $(f^*)^n=(f^n)^*$ on $H^{1,1}({  X})$.

\noindent{\it Proof.} It is sufficient to show that $(f^*)^2 = (f^2)^*$ on $Pic(X)$. If $D$ is a divisor, then $f^*D$ is the divisor on ${  X}$ which is the same as $f^{-1} D$ on ${  X} - {\cal I}$. Since ${\cal I}$ has codimension at least $2$, we also have $(f^2)^* D = f^*(f^* D)$ on ${  X} - {\cal I} - f^{-1} ({\cal I})$. By our hypothesis $f^{-1} ({\cal I})$ has codimension at least $2$. Thus we have $(f^2)^*D = (f^*)^2 D$ on ${  X}$. \qed

In a similar way, we may define $f^*:H^{p,q}({  X})\to H^{p,q}({  X})$.  That is, if  $\beta$ is a $(p,q)$ form on ${  X}$, then the pullback $\pi_2^*\beta$ is a smooth form on $\Gamma$.  We may  let $[\pi_2^*\beta]$ denote the reinterpretation of the form as a current, and we may push it forward to obtain a current $f^*\beta=\pi_{1*}[\pi_2^*\beta]$ on ${  X}$.  This pulls smooth forms back to currents and is well defined at the level of cohomology classes.  If $\alpha\in H^{p',q'}$ is an element of the dual cohomology group, then we have $\langle \alpha,f^*\beta\rangle=\langle\pi_1^*\alpha,\pi_2^*\beta\rangle$.  Now if $f$ is birational and $g=f^{-1}$, then  
$$\langle g^*\alpha,\beta\rangle=\langle\pi_1^*\alpha,\pi_2^*\beta\rangle=\langle\alpha,f^*\beta\rangle\eqno(1.1)$$  
If we have $(f^n)^*=(f^*)^n$ on $H^{p,q}$  for $n\ge0$ , then this gives us $(g^n)^*=(g^*)^n$ on $H^{p',q'}$.

The following is proved along standard lines: 
\proclaim Proposition 1.3.  If $f$ is a pseudo-automorphism, then we have $(f^n)^*=(f^*)^n$ on $H^{1,1}$ for all $n\in {\bf Z}$.  In particular, $(f^{-1})^*=(f^*)^{-1}$.

From this we get the following:

\proclaim Proposition 1.4. If $f:X\dasharrow X$ is a pseudo-automorphism on a d-dimensional manifold, then for all $n\in{\bf Z}$ we have $(f^*)^n=(f^n)^*$ on both $H^{1,1}$ and $H^{d-1,d-1}$.  Further, the characteristic polynomials of $f^*$ on $H^{1,1}$ and $H^{d-1,d-1}$ are the same and therefore the first and the $d-1^{\rm st}$ dynamical degree are the same.

\noindent {\it Proof.  } From Proposition 1.3, we have $f^* (f^{-1})^* = id$ and so we have that if $\alpha\in H^{1,1}$ and $\beta\in H^{d-1,d-1}$, then $f^*\alpha\cdot f^*\beta=((f^{-1})^*f^*)\alpha\cdot\beta=\alpha\cdot\beta$. Further, we conclude that $(f^*)^n=(f^n)^*$ on $H^{d-1,d-1}$ for all $n\in {\bf Z}$.
\qed

Now let us define some specific blowup situations.  This will serve to define the constructions we will use in the sequel, and it allows us to exhibit the models of indeterminate behavior that we will encounter.

\noindent{\it Blowing up a point and a line which contains it. }  We use $(x_0,x_1,x_2)\mapsto [x_0:x_1:x_2:1]$ as local coordinates in a neighborhood of $e_3:=[0:0:0:1]\in{\bf P}^3$.  Let ${  X}_1$ be the space obtained by blowing up a point $e_3$ and we let $E_3$ denote the fiber over $e_3$. We may use 
$$\pi_1:{  X}_1\ni (s_0, s_1, \xi_2)_1 \mapsto[\xi_2 s_0:\xi_2s_1:\xi_2:1] \in {\bf P}^3 \eqno{(1.2)}$$
as a local coordinate system for a neighborhood of $E_3 \cap \{ x_0=x_1=0\}$ in ${  X}_1$. It follows that the exceptional fiber $E_3 = \{ \xi_2=0\}$ in this coordinate system. 
%

%

Let $\Sigma_{01}=\{x_0=x_1=0\}\subset{\bf P}^3$ denote the $x_2$-axis.  The strict transform of $\Sigma_{01}$ inside ${  X}_1$ may be written as  $\Sigma_{01}=\{s_0=s_1=0\}$.  Thus $\Sigma_{01}\cap E_3=\{s_0=s_1=\xi_2=0\}$.  Let ${  X}_2$ be a complex manifold obtained by blowing up $\Sigma_{01}$ in ${  X}_1$. We can define a local coordinate system of ${  X}_2$ via $\pi_2 : {  X}_2 \ni (t_0, \eta_1, \xi_2)_2 \mapsto (t_0 \eta_1, \eta_1, \xi_2)_1 \in {  X}_1$
Thus $\pi_2\circ\pi_1:{  X}_2\to{\bf P}^3$ is given, in this coordinate neighborhood, by
$$\pi_1\circ\pi_2:{  X}_2\ni (t_0, \eta_1, \xi_2)_2\mapsto[t_0\eta_1\xi_2:\eta_1\xi_2:\xi_2:1] \in {\bf P}^3.\eqno(1.3)$$
The inverse of $\pi_1$ (resp.\ $\pi_2$) gives a model of indeterminate behavior that blows up the point $(0,0,0)$ (resp.\ the line $\{x_1=x_2=0\}$) to a hyperplane:
$$\pi_1^{-1}: (x_1,x_2,x_3)\mapsto (x_1/x_3,x_2/x_3,x_3), \ \ \ \pi_2^{-1}:(x_1,x_2,x_3)\mapsto (x_1/x_2,x_2,x_3). \eqno(1.4)$$

\noindent{\it Blowing up two intersecting lines. }  Let $\pi_1: {  Z}_1 \to {\bf P}^3$ be the blowup of the $x_1$-axis $\Sigma_{02}=\{x_0=x_2=0\} \subset {\bf P}^3$. We use local coordinate system in ${  Z}_1$ 
$$ \pi_1 : {  Z}_1 \ni (  \xi,x,s)_{{  Z}1} \mapsto [ s \xi: x: s:1] \in {\bf P}^3$$
Let us denote the blowup fiber over the point $o= \Sigma_{01} \cap \Sigma_{02}=[0:0:0:1] \in {\bf P}^3$ as ${\cal F}_o^1$ then in this coordinate system we have ${\cal F}_o^1 = \{ s=x=0\}$.  The strict transform of the $x_2$-axis in ${  Z}_1$ is given by $\ell_2= \{\xi=x=0\}$ and ${\cal F}_o^1 \cap \ell _2= ( 0,0,0)_{{  Z}_1}$. Now let ${  Z}_2$ be the blow up of $\ell_2$ with a local coordinate system $$\eqalign{\pi:= \pi_1\circ \pi_2: ( t, \eta, s)_{{  Z}_2} \in {  Z}_2 &\mapsto ( t, t\eta, s )_{{  Z}_1} \in {  Z}_1\cr & \mapsto [ ts:t \eta: s:1] \in {\bf P}^3}$$
We denote the second (new) fiber over $o$ as ${\cal F}_o^2$, so $ {\cal F}_o^2 =  ( 0,\eta,0)_{{  Z}_2} $. Let us also use ${\cal F}_o^1$ for its strict transform in ${  Z}_2$, so $ {\cal F}_o^1 \cup {\cal F}_o^2 = \pi^{-1}_2 \circ \pi_1^{-1} \{o\}$ and ${\cal F}^1_o =  ( t,0,0)_{{  Z}_2} $. 

Let $\tau[x_0:x_1:x_2:x_3]=[x_0:x_2:x_1:x_3]$ be the involution that interchanges the $x_1$- and $x_2$-axes.  It follows that $\tau$ induces the involution $\tilde\tau=\pi^{-1}\circ\tau\circ\pi$ on ${  Z}_2$.  In coordinates, we have
$$\tilde\tau: (t,\eta,s)\mapsto (s/\eta, \eta, t \eta), \eqno(1.5) $$
which will serve as our third model of indeterminate behavior.  We note that $\tilde\tau$ is regular on ${\cal F}^2_o-{\cal F}^1_o$, while each point of ${\cal F}^1_o$ blows up to the variety ${\cal F}^1_o$.

Similarly we can blowup the $x_2$- axis first and then the strict transform of $x_1$-axis. Performing similar computations, we obtain a blowup space $\hat\pi:{  Y}_2\to{\bf C}^3$.  The identity map $\iota$ on ${\bf P}^3$ lifts to a map $\tilde\iota:{  Z}_2\to{  Y}_2$, which in local coordinates is similar to $\tilde\tau$. 

\noindent

{\it Remark.}  Suppose that $\gamma'$ and $\gamma''$ are curves in ${\bf P}^3$ which intersect transversally at points $\{p_1,\dots,p_N\}$.   We have local coordinate systems for $1\le j\le N$ so that $p_j$ is the origin, and  $\gamma'$ (resp.\ $\gamma''$) coincides with the $x$-axis (resp.\ the $y$-axis) in a neighborhood of $p_j$.  Since the operation of blowing up the axes is local near $p_j$, we may construct a blowup space $\pi:W\to{\bf P}^3$ in which $\gamma'$ and $\gamma''$ are both blown up, and over each $p_j$ we are free to choose whether $\gamma'$ or $\gamma''$ was blown up first, independently of the choices over $p_k$ for $k\ne j$.

\proclaim  Theorem 1.5.   Let $f$ be a birational map of   $X$.  Let $X_0\subset X$ be a hypersurface such that the strict transform is $f(X_0)=X_0$.  Let $\varphi:X\to Y$ is a birational map which conjugates $(f,X)$ to an automorphism $(g,Y)$.  Then there is a birational map $\hat\varphi:X\to\hat Y$ such that the strict transform $\hat Y_0:=\hat\varphi( X_0)$ is a nonsingular hypersurface, and the induced map $\hat g:=\hat\varphi\circ f\circ\hat\varphi^{-1}$ gives an automorphsm of $\hat Y$.

\noindent{\it Proof.}  We may assume that $X_0$ is irreducible.  Since $X_0$ is a hypersurface, we may take its strict transform $\varphi(X_0)$.  If $\varphi(X_0)$ is a point in $Y$, then it is fixed by $g$.  If $\pi_1:Y_1\to Y$ be the blowup of the point $\varphi(X_0)$, then $g$ lifts to an automorphism of $Y_1$. Let $\phi_1:=\pi_1^{-1}\circ \varphi$.  If $\varphi_1(X_0)$ is again a point, we can repeat this blowing-up process until $\varphi_1(X_0)$ has dimension $>0$, which we may assume to be  1.  If the singular locus of $\varphi_1(X_0)$ is nonempty, it is finite and invariant under $f_1$.  Now we can blow up the singular set of $\varphi_1(X_0)$ finitely many times and have a new blow up space $\pi_2:Y_2\to Y_1$.  Since we were blowing up invariant sets, the induced birational map $g_2$ of $Y_2$ is again an automorphism.  Now the image $\varphi_2(X_0)$ must be a nonsingular curve, which must be invariant.  We can blow up this curve, and repeat the process finitely many times so that $\varphi_3(X_0)$ has dimension $>1$.  We continue this process until $\varphi_N(X_0)$ is a nonsingular hypersurface in $Y_N$, and now we set $\hat Y=Y_N$.  \qed

\proclaim Corollary 1.6.    Let $f$ be a birational map of   $X$.  Let $X_0\subset X$ be a hypersurface for which the strict transform is $f(X_0)=X_0$.  Let $\varphi:X\to Y$ is a birational map which conjugates $(f,X)$ to an automorphism $(g,Y)$. Then the restriction $(f_{X_0},X_0)$ is birationally conjugate to an automorphism.

\noindent {\it Proof of Theorem 4.}  Let $f$ be as in Theorem 4.  In \S {\bf C}  we study the restriction of $f^8$ to the plane $\Sigma_3=\{[x_0:x_1:x_2:x_3]\in {\bf P}^3: x_3=0\}$.  There we show that this restricted mapping is not birationally equivalent to an automorphism of $\Sigma_3$.    Thus Theorem 4 is a consequence of Corollary~1.6.  \qed

\noindent{\bf\S2.  Linear fractional recurrences.}  The maps (2.2) are among the Cremona transformations of 3-space which are discussed in Chapter 10 of [H].  We discuss general properties of these transformations, and for the generic parameters (2.10) we construct a new space $\pi:{  X}\to{\bf P}^3$, such that passing to the induced map $f_{  X}$ effectively eliminates one of the exceptional components.    

For $\{ i_1, \dots, i_k\} \subset \{0,1,2,3\}$,  we use the notation 
$$\Sigma_{i_1\cdots i_k} =\{x\in{\bf P}^3:x_{i_j}=0, 1\le j\le k\},\eqno(2.1)$$ and for a vector $A=(A_0,\dots, A_3)$ we will write $A\cdot x = A_0 x_0+A_1 x_1+A_2 x_2+A_3 x_3$. In homogeneous coordinates the maps (0.1) take the form 
$$f[x_0:x_1:x_2:x_3]=[x_0\beta\cdot x:x_2\beta\cdot x:x_3\beta\cdot x:x_0\alpha\cdot x]\eqno(2.2)$$ where $\alpha= (\alpha_0, \alpha_1,\alpha_2,\alpha_3)$, $\beta=(\beta_0,\beta_1,\beta_2,\beta_3)$. In the sequel, we will assume $$\alpha\ne\lambda\beta,\ \ \ \beta\ne(\beta_0,0,0,0), \ \ \ (\alpha_1,\beta_1)\ne(0,0).\eqno(2.3)$$ Note that if one of the first two conditions does not hold, then $f$ is linear, and if the third condition does not hold, then $f$ is independent of $x_1$ and thus $f$ is actually a 2-step recurrence.  If we set $\gamma=\beta_1\alpha-\alpha_1\beta$, then we have 
$${\cal I}=\Sigma_{\beta\gamma}\cup\Sigma_{0\beta}\cup\{e_1\}$$
where $\Sigma_{\beta}=\{\beta\cdot x=0\}$, $\Sigma_{\gamma} = \{\gamma\cdot x=0\}$, $\Sigma_0=\{x_0=0\}$, $\Sigma_{\beta\gamma}=\Sigma_\beta\cap\Sigma_\gamma$, $\Sigma_{0\beta}= \Sigma_0 \cap \Sigma_\beta$, and $e_1=[0:1:0:0]=\Sigma_{023}$.

The Jacobian determinant of $f$ is given by $2 x_0(\gamma \cdot x) (\beta \cdot x)^2$. Thus we see that the exceptional hypersurfaces are  ${\cal E} = \{ \Sigma_0, \Sigma_\beta, \Sigma_\gamma\}$. The action of $f$ on the exceptional varieties is given as follows: for $\lambda_2,\lambda_3 \in {\bf C}, (\lambda_2, \lambda_3) \ne (0,0)$,
$$f:\ \ \ \ \ \eqalign{&\Sigma_\beta  \mapsto e_3,\cr  &\Sigma_0 \cap \{ \lambda_2x_2= \lambda_3 x_3\} \mapsto [0: \lambda_3:\lambda_2:0] ,\cr &\Sigma_\gamma \cap \{ \lambda_2 x_2= \lambda_3 x_3\} \mapsto \Sigma_{BC} \cap \{ \lambda_2 x_1= \lambda_3 x_2\}\cr}\eqno{(2.4)}$$ where we set $\check \alpha = (\alpha_0, \alpha_2,\alpha_3,0)$, $\check\beta = (\beta_0,\beta_2,\beta_3,0)$,  and $$B=(-\alpha_1,0,0,\beta_1),\qquad C= \beta_1 \check \alpha-\alpha_1 \check \beta.$$ 
Thus $\Sigma_\beta$ is blown down to a point.  The pencil of lines in $\Sigma_\gamma$ passing through $e_1 \in \Sigma_0\cap \Sigma_\gamma$ are all mapped to points in $\Sigma_{BC}$. The pencil of lines in $\Sigma_0$ passing through $e_1$ are all mapped to points on the line $\Sigma_{03}$, which is again one of the exceptional lines.  We have strict transforms:
$$f: \Sigma_0 \mapsto \Sigma_{03} \mapsto e_1 \eqno{(2.5)}$$

The inverse is given by $$f^{-1}[x_0:x_1:x_2:x_3]=[x_0B\cdot x:x_0\check\alpha\cdot x-x_3\check\beta\cdot x:x_1B\cdot x:x_{2}B\cdot x],\eqno{(2.6)}$$
and the indeterminacy locus is ${\cal I}(f^{-1})=\Sigma_{0B}\cup\Sigma_{BC}\cup\{e_3\}.$ The Jacobian of $f^{-1}$ is $2x_0(C\cdot x)(B\cdot x)^{2}$, so the exceptional hypersurfaces are
${\cal E}(f^{-1})=\{\Sigma_0,\Sigma_B,\Sigma_C\}$ and for $\mu_1,\mu_2 \in {\bf C}, (\mu_1,\mu_2) \ne (0,0)$, $$f^{-1}:\ \ \ \ \eqalign{&\Sigma_B\mapsto e_1,\cr &\Sigma_0 \cap \{ \mu_1x_1=\mu_2 x_2\} \mapsto\Sigma_{0\beta} \cap \{ \mu_1x_2=\mu_2 x_3\},\cr &\Sigma_C\cap \{ \mu_1x_1=\mu_2 x_2\} \mapsto\Sigma_{\beta\gamma}\cap \{ \mu_1x_2=\mu_2 x_3\}.\cr}\eqno{(2.7)}$$
Now let us construct the space $\pi_1: {  X}_1 \to {\bf P}^3$ by blowing up a point $e_1$, and then  the space $\pi_2:{  X} \to {  X}_1$ obtained by blowing up a line $ \Sigma_{03}$. We set $$\pi= \pi_1\circ \pi_2 : {  X} \to {\bf P}^3. \eqno{(2.8)}$$ Let $S_{03}$ denote the blowup fiber over the strict transformation of $\Sigma_{03}$ in ${  X}_1$ and $E_1$ for the strict transformation of $ \pi_1^{-1} e_1$ in ${  X}_1$. For the induced map on ${  X}$, the orbit of $\Sigma_0$ becomes $$ f_{  X}: \Sigma_0 \to S_{03} \to E_1 \to \Sigma_B \eqno{(2.9)}$$ If $X$ and $Y$ are irreducible, we will say that a rational map $f:X \dasharrow Y$ is {\it dominant} if the rank of $df$ is equal to the dimension of $Y$ on a dense open set.  Let us define a generic condition:
 $$\beta_1 \ne 0,\ \ \beta_1 \alpha_2 \ne \alpha_1 \beta_2, \ {\rm \ and\ }\beta_1 \alpha_3 \ne \alpha_1 \beta_3.  \eqno(2.10)$$

For simplicity we use the same notation for both a variety and its strict transform, if there is no possibility of confusion. 

\proclaim Proposition 2.1. If (2.10) holds, then all the maps in $(2.9)$ are dominant, so ${\cal E}(f_{  X}) = \{ \Sigma_\beta, \Sigma_\gamma\}$.  Further, ${\cal  I} (f_{  X}) = \Sigma_{\beta 0}\cup \Sigma_{\beta\gamma}$.

\noindent {\it Proof.} Let us first consider the restriction to $S_{03}$. We may use the local coordinates for a neighborhood of $S_{03} $, $(s_0,x_2,\xi_3)_{S_{03}} \mapsto [s_0:1:x_2:s_0 \xi_3]$ . For the neighborhood of the exceptional fiber $E_1$ over $e_1$, we use $(t_0,\zeta_2, \zeta_3)_{E_1} \mapsto [t_0 : 1:t_0\zeta_2:t_0 \zeta_3]$. It follows that $S_{03} =\{ (0,x_2, \xi_3)_{S_{03}}\}$ and $E_1=\{ (0,\zeta_2,\zeta_3)_{E_1}\}$. Using these local coordinates we have
$$f_{  X}|_{S_{03}}: (0,x_2, \xi_3)_{S_{03}} \mapsto \left(0,\xi_3, {\alpha_1+ \alpha_2 x_2 \over \beta_1 +\beta_2 x_2} \right)_{E_1}. $$ 
To have a dominant map, it is required that $\beta_1 \alpha_2 \ne \alpha_1 \beta_2$. For the restrictions of the induced birational map to $\Sigma_0$ and $E_1$ are given by linear maps: 
$$\eqalign{ &f_{  X}: \Sigma_0 \ni [0:x_1:x_2:x_3] \mapsto \left(0, {x_3 \over x_2}, { \alpha_1 x_1+ \alpha_2 x_2 + \alpha_3 x_3 \over  \alpha_1 x_1+ \alpha_2 x_2 + \alpha_3 x_3} \right)_{S_{03}} \in S_{03}\cr &f_{  X}: E_1 \ni (0, \zeta_2, \zeta_3)_{E_1} \mapsto [\beta_1: \beta_1 \zeta_2: \beta_1 \zeta_3: \alpha_1] \in \Sigma_B.\cr}$$ We see that $f_{  X}|_{E_1}$ is dominant because $\beta_1 \ne 0$. And since  $\beta_1 \alpha_2 \ne \alpha_1 \beta_2$ and $\beta_1 \alpha_3 \ne \alpha_1 \beta_3$, we see that $f_{  X}|_{\Sigma_0}$ is dominant.\qed

Thus in passing to $f_{  X}$, we have removed one exceptional hypersurface and one point of indeterminacy.
There is a group of linear conjugacies acting on the family (0.2). For $(\lambda, c, \mu) \in {\bf C}^* \times {\bf C}^* \times {\bf C}$, we set  
$$(\alpha,\beta) \mapsto (\lambda\alpha, \lambda\beta) \eqno{(2.11a)}$$  
$$(\alpha,\beta) \mapsto (\alpha_0,c\alpha_1,c\alpha_2,c\alpha_3, c\beta_0, c^2\beta_1,c^2\beta_2,c^2\beta_3) \eqno{(2.11b)}$$ 
$$\eqalign{& \ \ \ \ \ \ \ \ \ \ (\alpha,\beta) \mapsto (\alpha_0',\alpha_1',\alpha_2',\alpha_3',\beta_0',\beta_1,\beta_2,\beta_3)\cr
\alpha_0'=  \alpha_0 &+\mu(\alpha_1+\alpha_2+\alpha_3) +\mu(\beta_0+\mu\beta_1+\mu\beta_2+\mu\beta_3),\ \  \alpha_1' = \alpha_1-\mu\beta_1,\cr 
  \alpha_2'=\alpha_2  &  -\mu\beta_2,\ \   \alpha_3'=\alpha_3-\mu\beta_3,\  \ \beta_0'= \beta_0+ \mu(\beta_1+\beta_2+\beta_3).\cr} \eqno{(2.11c)}$$ The first action corresponds to the homogeneity of $f$. The action (2.11b) corresponds to a scaling of $(x_1,x_2,x_3)$ in affine coordinates, and (2.11c) comes from translation by the vector $(\mu,\mu,\mu)$. Note that these actions preserve the form of the recurrence relation.

\medskip
\noindent{\bf\S3.  Non-critical maps.} A map $f$ of the form (2.2) is  {\it critical} if (3.1) holds:
$$\beta_2=\beta_3=0, {\rm\ \ and\ }\beta_1\alpha_2\alpha_3 \ne 0. \eqno(3.1)$$
In this section we establish the following:
\proclaim Theorem 3.1. If $f$ is not critical, then $f$ is not periodic. 

We will use the following criterion:
\proclaim  Proposition 3.2.  Suppose that $f:{  X}\to{  X}$ is periodic, i.e., $f^p_{  X}$ is the identity for some $p>1$.  If  $E\subset{  X}$ is an exceptional hypersurface, then $f^jE\subset{\cal I}(f_{  X})$ for some $1\le j<p$.

\noindent{\it Proof. }   Since $E$ is exceptional, then ${\rm codim}(f(E))\ge2$.  Let us consider the sequence of varieties $V_j:=f^j_{  X}(\Sigma)$.  If $V_j\not\subset{\cal I}(f_{  X})$ for all $j$, then applying the strict transform of $f$ repeatedly, we have $f_{  X}^{j+1}(E) = f_{  X}(V_j)$ for all $j$, so ${\rm codim}(f(V_j))\ge2$ for all $j$.  On the other hand, we must have $f^p_{  X}(E) = E = V_p$.  \qed

The proof of Theorem 3.1 will involve several cases, so we start with some Lemmas.
\proclaim Lemma 3.3.   Let $\pi:{  X} \to{\bf P}^3$ be the complex manifold defined in (2.8).  If $\beta_1 =0$, then there is a hypersurface $V \subset {  X}$ such that $f_{  X}^n V$ is a point of ${  X}-{\cal I}$ either for all $n\ge 1$ or for all $n \le -1$.  

\noindent {\it Proof.} We use the local coordinates for a neighborhood of $S_{03}$ and a neighborhood of $E_1$ defined in \S2. Since $\beta \ne (\beta_0,0,0,0)$, we first assume that $\beta_3 \ne 0$ (and thus we may suppose $\beta_3=1$) and to consider various cases.

{(i)} Case $\beta_2 \ne 0$ : In this case, the orbit of $\Sigma_\beta$ is:$$f_{  X} : \Sigma_\beta \mapsto e_3 \mapsto (0,0,\alpha_3)_{S_{03}} \mapsto (0,\alpha_3,  \alpha_2 / \beta_2  )_{E_1} \mapsto e_3 \in \Sigma_0 \setminus {\cal I}$$
Thus the orbit is periodic and remains a regular point of $f_X$.

{(ii)} Case $\beta_2 =0, \alpha_2 =0$. If $\beta_0 + \alpha_3 =0$, then we have $f_{  X}^2(e_3) = e_3$. If $\beta_0+\alpha_3 \ne 0$, then $f_X^3(e_3) = e_3$. In both cases, the orbit of $\Sigma_\beta $ is pre-periodic.

{(iii)} Case $\beta_2 =0, \alpha_2 \ne 0$. With this parameters, we have a two-cycle betweeen $\Sigma_{02}$ and the fiber over $e_2$, $\{ (0,0,\xi_3)_{S_{03}}: \xi_3 \in {\bf C}\}$.  $$f_{  X}: \ \ \ \eqalign{ &[0:x_1:0:x_3] \in\Sigma_{02} \mapsto (0,0, (\alpha_1 x_1 + \alpha_3 x_3)/x_3)_{S_{03}} \cr &(0,0,\xi_3)_{S_{03}} \mapsto  [0:\beta_0+\xi_3:0:\alpha_2] \in \Sigma_{02}.\cr}$$ Since $f^2_{  X}(e_3)=[0:\beta_0+\alpha_3:0:\alpha_2]$ and both $\Sigma_{02}$ and the fiber over $e_2$ in $S_{03}$ are disjoint from ${\cal I}$, the forward orbit of $\Sigma_\beta$ is disjoint from ${\cal I}$.

(iv)  Case $\beta_3 = 0$ and $\beta_2\ne 0$.  Under the backward mapping, $\Sigma_0$ is exceptional. To see the backward  iteration, let us use a different local coordinate system in a neighborhood of $E_1$: $(\zeta_0, \zeta_2, t_3)_{E_1} \mapsto [\zeta_0 t_3: 1: \zeta_2 t_3: t_3] \in {\bf P}^3$. Using this local coordinates we see 
$$ f^{-1}_{  X} : \eqalign{ &\Sigma_0 \ni [0: x_1:x_2:x_3] \mapsto  ( 0, x_1/x_2,0 )_{E_1} \in E_1\cap \Sigma_0 \cr
 & E_1 \ni ( \zeta_0, \zeta_2,0)_{E_1} \mapsto (0, \alpha_1 \zeta_0/(\beta_2 - \alpha_2 \zeta_0), \zeta_2 /\zeta_0)_{S_{03}}\in S_{03} \cr 
 & S_{03} \ni ( 0, x_2, \xi_3)_{S_{03}} \mapsto [ 0: -\alpha_2-\alpha_3 x_2 + \beta_2 \xi_3: \alpha_1: \alpha_1 x_2] \in \Sigma_0}$$
Let $p = \Sigma_0 \cap E_1 \cap S_{03}$. It follows that $f^{-1}_{  X} p =p$, that is, $p$ is a fixed point for the inverse mapping.  
\qed

Now let us suppose that $\beta_1 \ne0$. Using actions (2.11a--c), we may assume that $\beta_1=1$ and $\alpha_1=0$. 

\proclaim Lemma 3.4. Suppose that $\beta_1 =1, \alpha_1=0$. If either $\beta_2\ne0$ or $\beta_3 \ne 0$, then $\Sigma_0$ is exceptional and preperiodic for $f^{-1}$.  

\noindent {\it Proof.} Let us first assume that $\beta_3 \ne 0$. Then $$f^{-1}: \Sigma_0 \ni [0:x_1:x_2:x_3] \mapsto [0:-(\beta_2 x_1+\beta_3 x_2):x_1:x_2] \in \Sigma_{0\beta}$$ and $\Sigma_{0\beta}$ is invariant under $f^{-1}$ :
$$f^{-1}: [ 0:-(\beta_2 x_2+\beta_3 x_3):x_2:x_3] \mapsto [ 0:\beta_2 (\beta_2 x_2+\beta_3 x_3)-\beta_3 x_2:-(\beta_2 x_2+\beta_3 x_3):x_2].$$ Now suppose $\beta_3=0$ and $\beta_2 \ne 0$. In this case we have
$$f^{-1}: \Sigma_0 \ni [0:x_1:x_2:x_3] \mapsto [0:-\beta_2 x_1:x_1:x_2] \mapsto [0:\beta_2^2:-\beta_2:1]\in \Sigma_{0\beta}\setminus {\cal I}(f^{-1})$$
and this last point is fixed under $f^{-1}$. \qed

\medskip
Let $\pi:Z \to {\bf P}^3$ be the complex manifold obtained by blowing up $e_2$ and $\Sigma_{02}$ and let $E_2$ and $S_{02}$ be the corresponding blow-up divisors. In the following Lemma, we use the local coordinates $(s_0,x_1, \xi_2)_{S_{02}} \mapsto [s_0: x_1:s_0 \xi_2:1]$ in a neighborhood of $S_{02}=\{ s_0 =0\}$  and $(u_0, \eta_1,\eta_3)_{E_2} \mapsto [u_0:\eta_1 u_0:1:\eta_3 u_0]$ in a neighborhood of $E_2=\{u_0 =0\}$.
 
\proclaim Lemma 3.5. Suppose that $\beta_1 =1, \alpha_1=\beta_2=\beta_3=0$. If either $\alpha_2=0$ or $\alpha_3=0$ then $S_{02}$ is exceptional and pre-periodic for $f_Z$ or $f_Z^{-1}$.

\noindent{\it Proof.} If $\alpha_2=\alpha_3=0$ then the mapping is basically linear. \hfill \break
\item{(i)} Case $\alpha_3 = 0$ and $\alpha_2 \ne0$ : 
$$f_Z :  (0,t, \xi_2)_{S_{02}}\mapsto (0,\xi_2,0)_{E_2} \mapsto (0,\beta_0+\xi_2,0)_{S_{02}} \mapsto (0,0,0)_{E_2} \leftrightarrow (0,{\beta_0 \over \alpha_2},0)_{S_{02}}.$$
\item{(ii)} Case $\alpha_2 = 0$ and $\alpha_3 \ne0$ : 
$$f_Z^{-1}:(0,t, \xi_2)_{S_{02}}\mapsto (0,-\beta_0,\xi_2)_{E_2} \mapsto (0,-\alpha_3,-\beta_0)_{S_{02}} \leftrightarrow (0,-\beta_0,-\beta_0)_{E_2}.$$\qed

\proclaim Theorem 3.6. If $f$ is not critical, then there exists a complex manifold $X$ such that either there is an exceptional hypersurface $E \subset X$ for  an induced birational map $f_X$ such that $f^n_X E \not\subset {\cal I}(f) $ for $n=1,2,\dots\,$, or the analogous statement holds for $f_X^{-1}$.

\noindent{\it Proof.} Let $X$ denote either the space $X$ or $Z$ in the Lemmas above. This Theorem follows from the Lemmas 3.3--5. \qed

\noindent{\it Proof of Theorem 3.1.}  If  $f$ is not critical, then Theorem 3.6 says that in each case that there is an exceptional hypersurface that does not map into ${\cal I}(f_X)$.  By Proposition 3.2, then, $f$ is not periodic.   \qed

\noindent{\bf\S4.  Critical Maps.}  Here we study critical maps in general. Let us recall the condition for $f$ being critical : $\beta_2=\beta_3=0$ and $\beta_1\alpha_2\alpha_3 \ne0$. Using the action (2.11a--c) we may assume that a critical map satisfies:
 $$\beta_1 =1, \beta_2=\beta_3=0, \alpha_1=0, \alpha_2\ne0, \alpha_3=1. \eqno(4.1)$$ 
In this section, we show (Lemma 4.2) that for every critical map there is a blowup space $\pi:{  Y}\to {  X}$ such that the induced map $f_{  Y}$ has only one exceptional hypersurface, which is $\Sigma_\gamma$.   We determine the indeterminacy locus of $f_Y$ (Corollary 4.6) and the dynamical degree for the generic case (Theorem 4.8).
\proclaim Proposition 4.1. If $f$ is critical then $f^{-1}$ is conjugate to a critical map. 

\noindent{\it Proof.} Let $\beta= (\beta_0,1, 0,0)$ and $\alpha=(\alpha_0, 0,\alpha_2,1)$ be parameters of a critical map $f$. We consider a linear map $\phi:[x_0:x_1:x_2:x_3] \mapsto [x_0:x_3:x_2:x_1]$. It follows that we have $$\eqalign{ \phi^{-1} \circ f^{-1} \circ \phi :& [x_0:x_1:x_2:x_3]  \cr &\mapsto[x_0 x_1:x_2 x_1:x_3 x_1:x_0 ( \alpha_0 x_0 - \beta_0 x_1+x_2+\alpha_2 x_3)].\cr}$$
Thus $f^{-1}$ is conjugate to a critical map of the form (2.2) with parameter values  $\beta'=(0,1,0,0)$ and $\alpha'= ( \alpha_0, -\beta_0, 1,\alpha_2)$ which satisfy the condition $(3.1)$.  \qed

\noindent{\it Remark. } By Proposition 4.1, each result for $f$  corresponds to a result for $f^{-1}$.  The translation between $f$ and $f^{-1}$ is guided by notation: $\beta\leftrightarrow B$, $\gamma\leftrightarrow C$, $1\leftrightarrow3$: thus $f^{-1}\Sigma_C=\Sigma_{\beta\gamma}$,  etc.  
\medskip
%



If $(4.1)$ holds, it follows that  $e_3 = \Sigma_0 \cap \Sigma_\beta \cap \{ x_2=0\} \in {\cal I}$, and 
$$f: \Sigma_\beta\, \dasharrow \,e_3 \,\rightsquigarrow\, \Sigma_{01} \, \rightsquigarrow\, \Sigma_0 \,\dasharrow \,\Sigma_{03}\, \dasharrow \, e_1\, \rightsquigarrow\, \Sigma_B \eqno{(4.2)}$$  
We define a new complex manifold $\pi_{  Y}:{  Y}  \to {\bf P}^3$ by blowing up $e_1$ and $e_3$, then the strict transform of $ \Sigma_{01}$, followed by the strict transform of  $\Sigma_{03}$.  (Equivalently, we start with ${  X}$ and blow up the strict transform of $e_1$ and $\Sigma_{03}$.)  For $j=1,3$, we denote the exceptional divisor over $e_j$ by $E_j$ and the exceptional divisor over $\Sigma_{0j}$ by $S_{0j}$ for $j=1,3$. The induced birational map $f_{  Y}:{  Y} \to {  Y}$ maps 
$$f_{  Y} :  \Sigma_\beta \to E_3 \to S_{01} \to\Sigma_0 \to S_{03} \to E_1 \to \Sigma_B\eqno{(4.3)}$$

\proclaim Lemma 4.2. The maps in (4.3) are dominant; $\Sigma_\gamma$ is the unique exceptional hypersurface for $f_{  Y}$, and $\Sigma_C$ is the unique exceptional hypersurface for $f^{-1}_{  Y}$.

\noindent{\it Proof.} Using the local coordinates defined in \S2, we have
$$f_{  Y} :\ \ \  \eqalign{& \Sigma_\beta \setminus \Sigma_{\beta\gamma} \ni [x_0: -\beta_0 x_0 : x_2:x_3] \mapsto (0,{x_2 \over x_0}, {x_3 \over x_0})_{E_3} \in E_3\cr &E_3 \ni (0,\xi_1,\xi_2)_{E_3} \mapsto (0,\xi_2,\beta_0+\xi_1)_{S_{01}} \in S_{01}\cr & S_{01} \setminus\Sigma_{\beta\gamma}\ni (0,\eta_1,x_2)_{S_{01}} \mapsto [0:x_2 (\beta_0+\eta_1): (\beta_0+\eta_1):1+\alpha_2 x_2] \in \Sigma_0. \cr}$$
In Proposition 2.1 we showed that the maps $\Sigma_ 0\to S_{03} \to E_1 \to \Sigma_B$ are dominant. It follows that $\Sigma_\gamma$ is the only exceptional hypersurface for $f_{  Y}$, and $\Sigma_C$ is the only one for $f_{  Y}^{-1}$. \qed

For $p\in {\bf P}^3$, we will say that a  point of $\pi_{  Y}^{-1}p$ is at {\it level 1} if it could have been obtained by a blowup a point or curve in ${\bf P}^3$.  Thus the points of all fibers are of level 1, unless they lie over $e_1$, $e_3$, or $e_2=\Sigma_{01}\cap\Sigma_{03}$.  The fibers $E_1\cap S_{03}$ and $E_3\cap S_{01}$ represent the points of $E_1$ and $E_3$ which are at level 2.  Over $e_2$, we  define ${\cal F}_{e_2}^1:= S_{01}\cap \pi_{  Y}^{-1}e_2$ and ${\cal F}_{e_2}^2:=S_{03}\cap \pi_{  Y}^{-1}e_2$.  We see that ${\cal F}^i_{e_2}$, for $i=1,2$ is at level $i$.


%

We see that the three curves on level 2 are not indeterminate:
\proclaim Lemma 4.3. If $f$ is critical, then the indeterminacy loci ${\cal I}(f_{  Y})$ and ${\cal I}(f_{  Y}^{-1})$ do not contain $E_1\cap S_{03}$, $E_3\cap S_{01}$, or ${\cal F}^2_{e_2}$. 

\noindent{\it Proof.}   Let us first consider the blowup fiber over $E_3 \cap \Sigma_{01}$. For this fiber let us use a local coordinate $( \xi_0, \xi_1, t_2)_{E_3} \mapsto [t_2 \xi_0: t_2 \xi_1: t_2:1] \in {\bf P}^3$. It follows that  the strict transformation of $\Sigma_{01} = \{ (0,0,t_2)_{E_3} \}$ and $E_3 \cap \Sigma_{01} = (0,0,0)_{E_3}$. The local coordinates in a neighborhood of the second blowup fiber over $E_3 \cap S_{01}$  is given by  $(\eta_0, u_1 , t_2)_{E_3^{01}} \mapsto (\eta_0 u_1, u_1, t_2)_{E_3} \mapsto [ \eta_0 u_1 t_2: u_1 t_2: t_2:1] \in {\bf P}^3$ and we have the second blowup fiber $E_3 \cap S_{01} = \{ ( \eta_0, 0,0)_{E_3^{01}}\}$. With these coordinates, we see that 
$$ f_{  Y} : (\eta_0,0,0)_{E_3^{01}} \mapsto (0,0,\eta)_{S_{01}} = S_{01} \cap \Sigma_0$$
where $(\xi, t, x_3)_{S_{01}} \mapsto [ \xi t: t: 1: x_3]$ gives a local coordinates near $S_{01}$. It follows that the second blowup fiber $E_3 \cap S_{01}$ is not indeterminate for $f_{  Y}$. The computations for $f_{  Y}^{-1}$ and for $E_1\cap S_{03}$ are essentially the same, and we see that $E_1\cap S_{03}$ and $E_3\cap S_{01}$ are not indeterminate for  $f_{  Y}$ or to  $f_{  Y}^{-1}$.

To consider the second blowup fiber ${\cal F}^2_{e_2}$, we use local coordinates  $(\xi, s, x_3)_{01} \to [ \xi s: s: 1: x_3] $. In this coordinates we see that $S_{01} = \{ s=0\}$ and the strict transform of $\Sigma_{03} = \{ \xi =0, x_3=0\}$. Thus the local coordinates near the blowup of $\Sigma_{03}$ is given by  $( \eta, s, t) _{03} \mapsto ( \eta t, s, t)_{01} \mapsto [ \eta t s: s: 1: t]$ and we have ${\cal F}^2_{e_2} = \{ (\eta,0,0)_{03}\}$. With this coordinates, we have 
$$ f_{  Y} : {\cal F}^2_{e_2}  \ni(\eta_0,0,0)_{03} \mapsto (0,0,\alpha_2 \eta)_{E_1} = E_1 \cap \Sigma_0$$ 
where $(\xi_0, t_2, \xi_3)_{E_1} \mapsto [ \xi_0 t_2:1: t_2: \xi_3 t_2]$ is local coordinates near $E_1$. Similarly we see that $f_{  Y}^{-1}  {\cal F}^2_{e_2}  = E_3 \cap \Sigma_0$ and the mapping is dominant. \qed

Recall from \S2 that in ${\bf P}^3$ each point on $\Sigma_{\beta\gamma}$ blows up to a line in $\Sigma_C$.  Note that $[0:0:1:-\alpha_2] = \Sigma_{\beta\gamma}\cap \Sigma_{01}$, and let ${\cal F}_{0\beta\gamma}:=\pi_{  Y}^{-1} (\Sigma_{\beta\gamma}\cap \Sigma_{01})$.  Note that the base point is the intersection of $\Sigma_{01}$ and $\Sigma_{\beta\gamma}$, two indeterminate lines.  Similarly, we write ${\cal F}_{0BC}:=\pi_{  Y}^{-1}(\Sigma_{BC}\cap\Sigma_{03})$.


\proclaim Lemma 4.4. If $f$ is critical, the fiber ${\cal F} _{0\beta\gamma}$ is a component of ${\cal I}(f_{  Y})$, and we have $f_{  Y}\ :\ {\cal F}_{0\beta\gamma}\setminus \Sigma_{\beta\gamma} \rightsquigarrow {\cal F}_{0BC}.$

\noindent{\it Proof.} Let us consider a local coordinates in a neighborhood of the fiber ${\cal F}_{0\beta\gamma}$ and a local coordinates in a neighborhood of ${\cal F}_{0BC}$: 
$$\eqalign{ (s_0, \eta_1,x)_{S_{01}} &\sim [s_0:\eta_1 s_0:1: x ] \in {\bf P}^3 \ \ {\rm and\ \ } {\cal F}_{0\beta\gamma} = \{ s_0 = 0, x= -\alpha_2\} \cr (s_0,x, \eta_3)_{S_{03}} &\sim [s_0:1: x:\eta_3 s_0 ] \in {\bf P}^3 \ \ {\rm and\ \ } {\cal F}_{0BC} = \{ s_0 = 0, x= -\alpha_2\} .\cr }$$ 
Since $$f_{  Y} (s_0,\eta_1, x)_{S_{01}} = [s_0 (\beta_0 +\eta_1): \beta_0 +\eta_1: x (\beta_0 +\eta_1): \alpha_0 s_0 + \alpha_2 + x],$$
we see that for each $\eta_1$ we have 
$$\eqalign{f_{  Y}(0,\eta_1,-\alpha_2)_{S_{01}} =& \lim_{s_0 \to 0, x \to -\alpha_2} (s_0, x, {1 \over \beta_0 + \eta_1} { \alpha_0 s_0 + \alpha_2 + x \over s_0})_{S_{03}}\cr =& \{ (0, - \alpha_2, \eta)_{S_{03}}: \eta \in {\bf C}\}= {\cal F}_{0BC}}$$
\qed

Let us define a set $S\subset{  Y}$ to be {\it totally invariant} if  it is completely invariant for the total transform, or if for all $p\in S$, we have $({f_{  Y}})_*p\subset S$ and for all $p\notin S$ we have $({f_{  Y}})_*p\cap S=\emptyset$.

\proclaim Lemma 4.5.  If $f$ is critical, then $\Sigma_{02}$ is indeterminate for $f_{  Y}$.  Each point of $\Sigma_{02}$ blows up to ${\cal F}^1_{e_2}$, and ${\cal F}^1_{e_2}$ is mapped smoothly to $\Sigma_{02}$.  The set $\Sigma_{02}\cup{\cal F}^1_{e_2}$ is totally invariant.

\noindent{\it Proof. } Recall that $f \Sigma_{02} = e_2$, and the point $e_2$ was blown up.  We consider points $[s: 1: s \xi: x]$  which are close to $\Sigma_{02}$ when $s$ is small. We see that
$$ f_{  Y} : [s: 1: s \xi: x] \mapsto [ { s \over x} : { s \xi \over x} : 1: s { \alpha_0 s + \alpha_2 s \xi + x \over x(\beta _0 s +1)} ] = ( {1 \over \xi},{ s \xi \over x}, s { \alpha_0 s + \alpha_2 s \xi + x \over x(\beta _0 s +1)} )_{01}.$$ Letting $s \to 0$ we see that  $f_{  Y}  [0:1:0:x] \rightsquigarrow \{ ( \eta, 0,0)_{01} \}$. Using the same local coordinates we also see that
$$ f_{  Y} : {\cal F}^1_{e_2} \ni ( \eta, 0,0)_{01} \mapsto [0:1:0: { \alpha_2 \eta \over \beta_0 \eta +1} ] \in \Sigma_{02}.$$
For the second statement, we first notice that from (4.3) $f_{  Y} ( (\Sigma_{02}\cup{\cal F}^1_{e_2} )^c-  {\cal I}(f_{  Y}))$ is disjoint from the set $\Sigma_{02}\cup{\cal F}^1_{e_2} $.  Since ${\cal I}(f_{  Y}) = \Sigma_{\beta\gamma} \cup {\cal F}_{0\beta\gamma} \cup \Sigma_{02}$ and $\Sigma_{\beta\gamma} \cap \Sigma_{02} = \emptyset$, we see that the set $\Sigma_{02}\cup{\cal F}^1_{e_2}$ and $\Sigma_{\beta\gamma} \cup {\cal F}_{0\beta\gamma}$ are disjoint. It follows that the set $\Sigma_{02}\cup{\cal F}^1_{e_2}$ is totally invariant  
\qed

\proclaim Corollary 4.6.  If $f$ is critical, then ${\cal I}(f_{  Y})=\Sigma_{\beta \gamma}\cup {\cal F}_{0\beta\gamma}\cup \Sigma_{02}$ has pure dimension 1.  

The behavior of $f_{  Y}$ at $\Sigma_{02}$ is, in suitable coordinates, given by the third model (1.5).   The behavior of $f_{  Y}$ at ${\cal F}_{0\beta\gamma}$, as seen in Lemma 4.4,  is different from the model (1.5).  Further, we note that by Proposition 4.1 and the remark following it, the analogues of Lemmas 4.2--5 all hold for $f^{-1}_{  Y}$.  For instance,  $\Sigma_C$ is the unique exceptional hypersurface for $f^{-1}_{  Y}$,  ${\cal I}(f^{-1}_{  Y})={\cal F}^1_{e_2}\cup {\cal F}_{0BC}\cup \Sigma_{BC}$, and each point of ${\cal F}_{0BC}-\Sigma_{BC}$ blows up under $f^{-1}_{  Y}$ to ${\cal F}_{0\beta\gamma}$.

\proclaim Corollary 4.7.  If $f$ is critical, then $f^j_{  Y}\Sigma_\gamma\cap(\Sigma_{02}\cup{\cal F}^1_{e_2})=\emptyset$ for all $j\ge0$.

\noindent{\it Proof. }  By Lemma 4.5, it suffices to consider the case $j=0$.  By $(4.1)$, $e_2\notin\Sigma_\gamma$ in ${\bf P}^3$, so the fiber over $e_2$ remains disjoint from $\Sigma_\gamma$ inside ${  Y}$.  Now $e_1=\Sigma_{02}\cap\Sigma_\gamma$ in ${\bf P}^3$ and we see that $\Sigma_{02}$ and $\Sigma_\gamma$ are separated when we blow up $e_1$ to make ${  Y}$. \qed

Recall that the degree complexity  is $\delta(f)=\lim_{n\to\infty}({\rm deg}(f^n))^{1/n}$.   If $\delta(f)>1$, then the degrees of the iterates $f^n$ grow exponentially in $n$.  In particular, $f$ cannot be periodic if $\delta(f)>1$.
  
\proclaim Theorem 4.8. If $f$ is critical, and if $f^n_{  Y} \Sigma_\gamma\not\subset \Sigma_{\beta\gamma} \cup {\cal F}_{0\beta\gamma}$,  then the first dynamical degree is $\delta(f) \sim 1.32472$, the largest root of $x^3-x-1$.

\noindent{\it Proof.}  Using Corollary 4.7 we see that $f^m_{  Y} \Sigma_\gamma \cap {\cal I}(f_{  Y}) = \emptyset$ for all $m \ge 1$. Thus by Proposition 1.2 we have  $(f_{  Y}^*)^m = (f_{  Y}^m)^*$ for all $m$. Thus $\delta(f)$ is the spectral radius of $f_{  Y}^*$. Inside the Picard group $Pic({  Y})$, we let $H_{  Y}$ be the class of  a generic hyperplane in ${  Y}$, and we  have
$$f^*_{  Y} :\  \eqalign{& H_{  Y} \to 2 H_{  Y} - E_1-E_3-S_{01} ,\  \ S_{01} \to E_3 \to \Sigma_\beta = H_{  Y} - E_3-S_{01},\cr 
  &E_1 \to S_{03} \to\Sigma_0 = H_{  Y} -E_1-E_3 -S_{01} - S_{03} \cr}  \eqno(4.4)$$ 
The characteristic polynomial of this transformation is $(x^2+1) (x^3-x-1)$, so $\delta(f)$ is as claimed.  \qed

Now we give the existence of Green currents, which are invariant currents with the equidistribution properties given in the following: 

\proclaim Theorem 4.9.   If $f$ is as in Theorem 4.8, then there is a positive closed current $T^+_Y$ in the class of $\alpha^+_Y$ with the properties: $f^*_YT^+_Y=T^+_Y$, and if \ $\Xi^+$\/ is a smooth form which represents $\alpha^+_Y$, then $\lim_{n\to\infty}\delta_1(f)^{-n} f_{Y}^{ n*}\Xi^+_Y =T^+_Y$  in the weak sense of currents on $Y$.


\noindent{\it Proof.}  Recall from Corollary 4.6 that ${\cal I}(f_Y)=\Sigma_{02}\cup \Sigma_{\beta\gamma}\cup{\cal F}_{0\beta\gamma}$.  The total forward image of this set is $\pi_2\pi_1^{-1}{\cal I}(f_Y)={\cal F}^1_{e_2}\cup {\cal F}'_{0BC}\cup\Sigma_C$.  We will show that if $\sigma\subset\pi_2\pi_1^{-1}{\cal I}(f_Y)$ is any curve, then $\alpha^+_Y\cdot\sigma\ge0$.  The Theorem will then be a consequence of Theorem 1.3 of [Ba].  

Up to a scalar multiple, we may write $\alpha=H_Y - c_1 E_1 - c_3 E_3 - c_{01}S_{01} - c_{03} S_{03}$.  Then since $f_Y^*$ is given by (4.4) we have $1>c_1>c_3> c_{01}=c_{03}>0$, and $c_1+c_3=1$.   Let us start with ${\cal F}_{0\beta\gamma}\subset{\cal I}(f_Y)$.  Points of this curve are blown up to ${\cal F}_{0BC}$.  The curve $\sigma={\cal F}_{0BC}$ is the exceptional fiber inside $S_{03}$ over the point $\Sigma_{BC}\cap \Sigma_{03}\in {\bf P}^3$.  Thus  $\sigma\cdot S_{03}=-1$, so $\alpha\cdot\sigma=c_{03}>0$.   Points of the indeterminate curve $\Sigma_{02}$ blow up to $\sigma={\cal F}^1_{e_2}$.  In this case, we have that $\sigma\cdot S_{01} $ and $\sigma\cdot S_{03}$, are $\pm1$, with opposite signs, depending on the order of blow-up of $\Sigma_{01}$ and $\Sigma_{03}$.  Thus $\sigma\cdot\alpha^+_Y = \pm c_{01}  \mp c_{03} = 0$

The other possibility is that  $\sigma\subset\Sigma_C$.  In this case, we have $\sigma\cdot H={\rm deg}(\sigma)= \sigma\cdot S_{01}$.  Further, if we let $m_3$ denote the multiplicity of $\sigma$ at $e_3$, then $m_3$ is bounded above by ${\rm deg}(\sigma)$.  If $\sigma$ is a curve in $\Sigma_C$, then it is represented by $L + m_1{\cal F}^1_{01} + m_3\epsilon_3$, where ${\cal F}^1_{01}$ represents a fiber of $S_{01}$, and $\epsilon_3=E_3\cap\Sigma_C$.  The multiplicities  $m_1$ and $m_3$ are bounded below by $-{\rm deg}(\sigma)$.  Since ${\cal F}^1_{01}\cdot S_{01}=-1= E_3\cdot \epsilon_3$, we have $\sigma\cdot \alpha^+_Y\ge {\rm deg}(\sigma)(1 - c_{01}-c_3)>0$.
\qed

\proclaim Lemma 4.10.  Suppose that $\Sigma_\gamma$ is not preperiodic.   If either $f^j_{  Y}\Sigma_\gamma={\cal F}_{0\beta\gamma}$ or $f^j_{  Y}\Sigma_\gamma$ is a point of ${\cal F}_{0\beta\gamma}-\Sigma_{\beta\gamma}$, then $f^{j+1}_{  Y}|_{\Sigma_\gamma}$ has rank 1, and $f^{j+1}_{  Y}\Sigma_\gamma={\cal F}_{0BC}$.

\noindent{\it Proof. } Let us describe the indeterminate behavior of $f_{  Y}$ at $\Sigma_{\beta\gamma}$ and ${\cal F}_{0\beta\gamma}$.  Up to coordinate changes in domain and range, we may assume that the indeterminate curve is $\{\xi_2=\xi_3=0\}$, and the maps behave like
$$\eqalign{& (s_1,\xi_2,\xi_3) \mapsto (s_1, {\xi_2 \over \xi_3}, \xi_3) \ \ \ {\rm near\ } \Sigma_{\beta\gamma}  \cr 
& (s_1, \xi_2, \xi_3) \mapsto ( s_1{ \xi_2 \over \xi_3}, \xi_2, \xi_3)\ \ \ {\rm near\ } {\cal F}_{0\beta\gamma}-\Sigma_{\beta\gamma}.}$$ 
The behavior near $\Sigma_{\beta\gamma}$ is given in (2.7), and ${\cal F}_{0\beta\gamma}$ is given in Lemma 4.4.  

We will track the forward orbit $f^i_{  Y}\Sigma_\gamma$.  Without loss of generality, we may assume that $f^i_{  Y}\Sigma_\gamma$ is not $\Sigma_{\beta\gamma}$ for $1\le i\le j$.  Let us choose coordinates $(u_1,u_2,v)$ such that $\Sigma_\gamma=\{v=0\}$, and the exceptional fibers of $\Sigma_\gamma$ are $\{v=0,u_1=const\}$.  Thus $f$ has the form $( u_1, u_2, v) \mapsto ( u_1, u_2 v,v)$.  By making coordinate changes in the range, we may represent the iterated map near $\Sigma_\gamma$ as $(u_1, u_2 v, v)$ as long as $f_{  Y}^i\Sigma_\gamma$ is not an exceptional curve in $\Sigma_\gamma$, and not a component of ${\cal I}(f_{  Y})$. 

Now suppose that $f^i_{  Y}\Sigma_\gamma=f^{i-1}_{  Y}\Sigma_{BC}$ coincides with ${\cal F}_{0\beta\gamma}$, which we will write as $\{\xi_2=\xi_3=0\}$, as above.  Thus we may assume that $(s,\xi_2,\xi_3)=(u_1,u_2 v,v)$, so $f^{i+1}$ has the form $(u_1,u_2,v)\mapsto (u_1u_2,u_2v,v)$, which is a map of rank 1.

As we continue to iterate $f_{  Y}$,  the other possibility is that $f^i_{  Y}\Sigma_\gamma\subset\Sigma_\gamma$ is an exceptional curve.  The coordinate of the map which varies when $v=0$ must be inside $\{u_1=const\}\subset\Sigma_\gamma$, which means that the map must be like $\Sigma_\gamma\ni (u_2v,u_1u_2,v)\mapsto (u_2v,u_1u_2v,v)$, which belongs to $\Sigma_{BC}$.  Now we continue to iterate this point forward.  It cannot re-enter $\Sigma_\gamma$, because otherwise the orbit would re-enter $\Sigma_{BC}$ and become pre-periodic.  Thus the only possibility is to enter the indeterminacy locus.  Let $N$ denote the first positive integer for which $f^N_{  Y}\Sigma_\gamma\subset\Sigma_\gamma$ was an exceptional curve. Then ${\cal I}(f_{  Y})\cap \bigcup_{i=0}^{N-1}\Sigma_{BC}$ is finite.  Thus the forward orbit of the point can intersect the indeterminacy locus only finitely many times.  If the point enters $\Sigma_{\beta\gamma}=\{\xi_2=\xi_3=0\}$, then the next image is  $(s,\xi_2,\xi_3)=(u_2v,u_1u_2v,v)\mapsto (u_2v,u_1u_2,v)$, which has rank 1 again, and we continue as before. \qed



\bigskip
\noindent{\bf\S5.  Pseudo-automorphisms.} In this section we assume that $f$ is critical and we consider the condition
$$f^N_{  Y}\Sigma_\gamma = \Sigma_{\beta\gamma} \ \ {\rm for\ some\ } N.  \eqno(5.1) $$
We give conditions for $f_{  Y}$ to be birationally equivalent to a pseudo-automorphism  (Theorem 5.1). 

Suppose that $(5.1)$ holds. For $1 \le j \le N-1$, we consider four possibilities: 
\itemitem{(i)} $f_{  Y}^j \Sigma_\gamma $ is a one of the exceptional fibers in $\Sigma_\gamma$, that is $f_{  Y}^j \Sigma_\gamma = \Sigma_\gamma \cap \{ \lambda_2 x_2 = \lambda_3 x_3\}$ for some $\lambda_2, \lambda_3 \in {\bf C}$ and $f_{  Y}^{j+1} \Sigma_\gamma$ is a point in $\Sigma_{BC}$. 
\itemitem{(ii)} $f_{  Y}^j \Sigma_\gamma $ is a point of indeterminacy in $\Sigma_{\beta\gamma}$, and this point blows up to one of fibers of $\Sigma_C$.
\itemitem{(iii)} $f_{  Y}^j \Sigma_\gamma \subset{\cal F}_{0\beta\gamma} $, and $f_{  Y}^{j+1} \Sigma_\gamma = {\cal F}_{0BC}$. 
\itemitem{(iv)} None of the above. 


\proclaim Theorem 5.1.  Suppose that a critical map $f$ satisfies $(5.1)$, and that whenever  case (iii) occurs above, then $f_{  Y}^j\Sigma_\gamma={\cal F}_{0\beta\gamma}$.  Then there is a blowup space $\pi:{  Z}\to{  Y}$ such that $f_{  Z}$ is a pseudo-automorphism. 

\noindent{\it Proof.} Let us define two sets of subvarieties:  $\Lambda_1 = \{ f_{  Y}^j \Sigma_\gamma\ |\  {\rm dim} f_{  Y}^j \Sigma_\gamma = 0, 1 \le j \le N\}$ and $\Lambda_2 = \{ f_{  Y}^j \Sigma_\gamma\ |\  {\rm dim} f_{  Y}^j \Sigma_\gamma = 1, 1 \le j \le N\}$. We construct a blowup space $\pi_{  Z}:{  Z}\to {  Y}$ obtained by first blowing up points in $\Lambda_1$ and then blowing up curves in $\Lambda_2$. We denote by ${\cal F}_1 , \dots, {\cal F}_N$ the blowup fibers over $f_{  Y} \Sigma_\gamma= \Sigma_{BC}, f_{  Y}^2 \Sigma_\gamma, \dots, f_{  Y}^N \Sigma_\gamma= \Sigma_{\beta\gamma}$, respectively.

We claim that the induced map $f_{  Z}$ on the orbit of $\Sigma_\gamma$ are dominant. To show this we need to check the map $f_{  Z}$ on $S_j$ where $f^j_{  Y}\Sigma_\gamma \subset {\cal E}(f_{  Y})$ or $f^j_{  Y}\Sigma_\gamma \subset {\cal I}(f_{  Y})$. For ${\cal F}_1$ we use local coordinates 
$$(x_1, \xi, s)_1 \mapsto [ 1:x_1: -\alpha_0 -\alpha_2 x_1 + s \xi: s] \in {\bf P}^3, \quad {\rm and \ \ }\{ s=0\} = {\cal F}_1.$$ Using this  we see that $$f_{  Z} : \Sigma_\gamma \ni [ x_0:x_1:x_2: - \alpha_0 x_0 - \alpha_2 x_2] \mapsto ( {x_2 \over x_0}\ ,\  \beta_0 + {x_1 \over x_0}\ ,\ 0)_1 \in {\cal F}_1.\eqno{(5.2)}$$ Suppose the possibility (i) occurs, that is $f_{  Y}^j \Sigma_\gamma=\Sigma_\gamma\cap \{ x_2 = \lambda x_3\}$ for some $\lambda \in {\bf C}$. Local coordinates near ${\cal F}_j$ and ${\cal F}_{j+1}$ are given by
$$ ( s, x, \xi)_j \mapsto [ -(\alpha_2 \lambda+1)/\alpha_0 + s : x : \lambda + s \xi: 1] \in {\bf P}^3, \qquad {\rm and\ \ }\{s=0\} = {\cal F}_j$$
$$ (s, \xi_1, \xi_2)_{j+1} \mapsto [  -(\alpha_2 \lambda+1)/\alpha_0 + s : \lambda + s \xi_1: 1 : s \xi_2]  \in {\bf P}^3, \quad {\rm and\ \ }\{s=0\} = {\cal F}_{j+1}.$$ It follows that 
$$ f_{  Z} : {\cal F}_j \ni ( 0, x, \xi)_j \mapsto (0\ ,\ \xi\ ,\  {(1+ \alpha_2 \lambda) ( \alpha_0+ \alpha_2 \xi)\over \beta_0 -\alpha_0 x + \alpha_2 \beta_0 \lambda} ) _{j+1} \in {\cal F}_{j+1}.\eqno{(5.3)}$$ In case the second possibility (ii), i.e. $f_{  Y}^j \Sigma_\gamma$ is a point in $\Sigma_{\beta\gamma}$ occurs is essentially identical to the first possibility (i). For the third possibility (iii), due to Lemma 6.3 we only need to check  the induced map on the blowup fiber over $f_{  Y}^j \Sigma_\gamma = {\cal F}_{0\beta\gamma}$. We use local coordinates near ${\cal F}_j$, the blowup fiber over ${\cal F}_{0\beta\gamma}$, and ${\cal F}_{j+1}$, the blowup fiber over ${\cal F}_{0BC}$:
$$ (s, \eta_1, \eta_2)_j \mapsto (s, \eta_1+s \eta_2, -1/\alpha_2)_{S_{0,1}} \mapsto [s: s\eta_1+ s^2 \eta_2: -1/\alpha_2:1]\in{\bf P}^3$$ 
$$ (s, \eta_1, \eta_2)_{j+1} \mapsto (s, -\alpha_2,\eta_1+s \eta_2)_{S_{0,3}} \mapsto [s:1:-\alpha_2: s\eta_1+ s^2 \eta_2]\in{\bf P}^3.$$ Using these local coordinates we see that 
$$f_{  Z} : {\cal F}_j \ni (0,\eta_1, \eta_2)_j \mapsto (0,{\alpha_0 \over \beta_0 + \eta_1}, { \alpha_0 \eta_2 \over \alpha_2 (\beta_0 + \eta_1)^2})_{j+1} \in {\cal F}_{j+1}.\eqno{(5.4)}$$ Now for the case when $f_{  Y}^j \Sigma_\gamma= E_3 \cap \Sigma_0$, $f_{  Y}^{j+1} \Sigma_\gamma= E_2 \cap \Sigma_0$,$f_{  Y}^{j+2} \Sigma_\gamma= E_1 \cap \Sigma_0$ we use local coordinates
$$ (\xi_0, \xi_1, t)_j \mapsto ( t \xi_0, \xi_1, t)_{E_3} \mapsto [ t^2 \xi_0: t \xi_1 : t:1] \in {\bf P}^3 $$ 
$$ (\xi_0, \xi_1, t)_j \mapsto ( t \xi_0, \xi_1, t)_{E_2} \mapsto [ t^2 \xi_0: t \xi_1 : 1:t] \in {\bf P}^3 $$ 
$$ (\xi_0, \xi_1, t)_j \mapsto ( t \xi_0, \xi_1, t)_{E_1} \mapsto [ t^2 \xi_0: 1:t \xi_1 : t] \in {\bf P}^3. $$
We have 
$$f_{  Z}  : {\cal F}_j \ni ( \xi_0, \xi_1, 0)_j \mapsto \left( {\xi_1^2 \over \xi_0}, { \xi_1 \over \xi_0},0  \right)_{j+1} \in {\cal F}_{j+1}\eqno{(5.5)}$$
 $$f_{  Z}  : {\cal F}_{j+1} \ni ( \xi_0, \xi_1, 0)_j \mapsto \left ( {\xi_1^2 \over\alpha_2^2 \xi_0}, { \xi_1 \over \alpha_2\xi_0},0 \right)_{j+1} \in {\cal F}_{j+1}.\eqno{(5.6)}$$
 The last part we have to check is $f_{  Z}$ on ${\cal F}_N$, the blowup fiber of the line of indeterminacy $\Sigma_{\beta\gamma}$. The local coordinates we use near ${\cal F}_N$ is given by
 $$ ( s, x_2, \xi_3)_N \mapsto [ 1: -\beta_0 + s: x_2 : -\alpha_0 - \alpha_2 x_2 + s \xi_3] \in {\bf P}^3$$ and we get
 $$f_{  Z} : {\cal F}_N \ni (0, x_2, \xi_3)_N \mapsto [ 1: x_2: -\alpha_0-\alpha_2 x_2: \xi_3] \in \Sigma_C.\eqno{(5.7)}$$ From (5.2--7) we see that the induced mapping $f_{  Z}$ is dominant on the orbit of $\Sigma_\gamma$ and therefore $f_Z$ has no exceptional hypersurface.  By Lemma 4.1, it follows that $f^{-1}_{  Z}$ also has no exceptional hypersurface, so $f_{  Z}$ is a pseudo-automorphism.  \qed

\proclaim Lemma 5.2.  Suppose that a critical map $f$ satisfies $(5.1)$, and that whenever  case (iii) occurs above, then the possibility $(iii)$ can occur at most once. 

\noindent{\it Proof.} Suppose that there are $0<j_1<j_2<N$ such that $f_{  Y}^{j_1} \Sigma_\gamma, f_{  Y}^{j_2} \Sigma_\gamma \subset {\cal F}_{0\beta\gamma}$. It follows that $f_{  Y}^{j_1+1}\Sigma_\gamma = f_{  Y}^{j_2+1}\Sigma_\gamma = {\cal F}_{0BC}$ and thus ${\cal F}_{0BC}$ is fixed under $f_{  Y}^{j_2-j_1}$. Since $0<j_1<j_2<N$, we have that the dimension $f_{  Y}^{j} \Sigma_\gamma \le1$ for all $j= j_1+1, \dots, j_2.$ Using the fact that ${\cal F}_{0BC}$ is fixed under $f_{  Y}^{j_2-j_1}$, we conclude that the dimension $f_{  Y}^j \Sigma_\gamma \le 1$ for all $j\ge 1$ which contradicts to the assumption that $f$ satisfies $(5.1)$.

\qed 

\medskip
\epsfysize=0.65in
\centerline{ \epsfbox{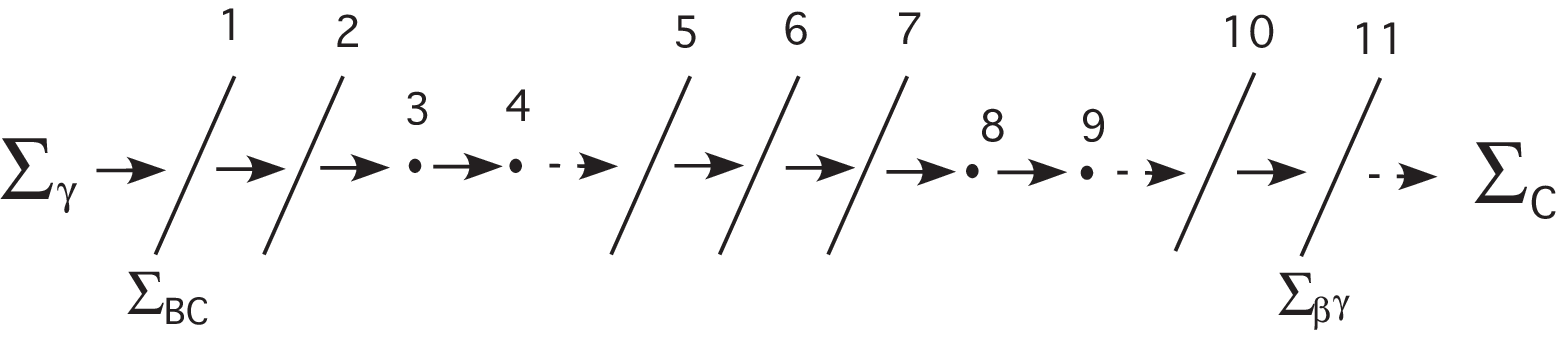} }

\smallskip
\centerline{Figure 5.1.  A hypothetical orbit: $d_1=2, u_1=4, d_2=7, u_2=9, m_d=m_u=2, N=11$.}
\medskip

Let $N$ be the smallest positive integer such that $f_{  Y}^N\Sigma_\gamma = \Sigma_{\beta\gamma}$. Let $m_d$ positive integers $d_1< d_2< \cdots < d_{m_d}$ denote the number of iteration to have the possibility (i), that is for each $j=1, \dots, m_d$,  $f_{  Y}^{d_j} \Sigma_\gamma$ is $\Sigma_\gamma \cap \{ \lambda_2 x_2 = \lambda_3 x_3\}$ for some $\lambda_2, \lambda_3 \in {\bf C}$. Let $u_1 < u_2< \cdots< u_{m_u}$ be $m_u$ positive integers such that $f_{  Y}^{u_j}\Sigma_\beta$ is of possibility (ii) for $j = 1, \dots, m_u$. We also set $m_s$ be a positive integer such that $f_{  Y}^{m_s +2} \Sigma_\gamma \subset {\cal  F}_{0\beta\gamma}$
if such possibility occurs. If there is no such case we set $m_s= \infty$.  To illustrate this numbering scheme, a hypothetical orbit of $\Sigma_\gamma$ is given in Figure 5.1. Here we have assumed that we are in the simpler case $m_s=\infty$, which means that the orbit never enters ${\cal F}_{0\beta\gamma}$, so the case (iii) does not occur.  Thus in Figure 5.1 the dimension can increase from 0 to 1 only via case (ii).   

Let us use the numbers $m_s$,  $m_u$, $m_d$, $u_j$, $d_j$ and $N$ to define four Laurent polynomials:
$$\eqalign{  Q_1  & := -1 - \sum_{j=1}^{m_d} {1 \over t^{d_j}}+ {1 \over t^{ m_s+1}}  \cr
 Q_2  & := {1 \over t^{m_s}} \left( {1 \over t} + {1 \over t^2} + { 1\over t^3} + { 1 \over t^4} \right)  \cr
Q_3  &  := -1 - \sum_{j=1}^{m_d} {1 \over t^{d_j}}+ {1 \over t^{ m_s}} \left( 1 + {1 \over t} - { 1 \over t^4}\right)  \cr 
Q_4 & := -t -t \sum_{j=1}^{m_d} {1 \over t^{d_j}} - t \sum_{j=1}^{m_u} {1 \over t^{u_j}}-{ 1\over t^{N-1}}-{1 \over t^{m_s}} \left( {1 \over t} + { 1\over t^3} \right).}$$

\proclaim Theorem 5.3. If  $f$ is pseudo-automorphism and if the possibility $(iii)$ can occur at most once, then the dynamical degree of $f$ is given by the largest root of the polynomial  
$$  \chi_f(t): = t^{N-1}\left[ (Q_1 - Q_4) t^3  +( 2 Q_1 - Q_2 - Q_3 -Q_4) t^2 + ( Q_1 - Q_3) t + Q_4 \right].\eqno{(5.8)}$$

\noindent{\it Proof.}  By Corollary 6.2 and Lemma 6.3, we see that $f_{  Y}$ satisfies the hypotheses of Theorem 5.1.  Let  $f_{  Z}$ be the corresponding pseudo-automorphism.  The dynamical degree will be the modulus of the largest root of the characteristic polynomial of $f_{  Z}^*$.  In the Appendix we show that the characteristic polynomialis given by $\chi_f$.  \qed


%
%

\bigskip
\noindent{\bf\S6.  Periodic maps.} In this section, we determine all possible periodic 3-step recurrences.  By \S3, we may assume $(4.1)$.  The question of periodicities for maps $(4.1)$ with $\beta_0=0$ has been answered by Cs\"ornyei and Laczkovic [CL]: they have shown that the only periodicities in this case are the two period 8 maps given in the Theorem stated in the Introduction.  We will consider the general case where $\beta_0$ is possibly nonzero.  We start by giving a necessary condition for a map to be periodic.

\proclaim Proposition 6.1. If $f$ is pseudo-automorphism and if $E$ is an exceptional hypersurface then there is an exceptional hypersurface $E'$ for $f^{-1}$ such that $f^n E= E'$ for some $n>0$ and the co-dimension of $f^jE$ is $\ge2$ for all $j=1, \dots, n-1$. 

\noindent{\it Proof.} Suppose $f$ has period $p$. Since $f^pE=E$ and ${\rm codim\ }fE \ge2$, it follows that there exists $0 <n\le p$ such that ${\rm codim\ }f^{n-1}E \ge2$ and ${\rm codim\ }f^nE =1$. Thus $f^nE$ is an exceptional for $f^{-1}$.\qed

Since
$f$ is critical,  ${\rm dim} f^j \Sigma_\beta <2$ for $j=1,2$, and  $f^3 \Sigma_\beta = \Sigma_0$; further,  ${\rm dim} f^j \Sigma_0 <2$ for $j=1,2$, and  $f^3 \Sigma_0 = \Sigma_B$.  By Lemma 4.2  the only exceptional hypersurface for $f_{  Y}$ is $\Sigma_\gamma$, and the only exceptional hypersurface for $f^{-1}_{  Y}$ is $\Sigma_C$.  This gives us the following necessary condition for $f$ to be periodic. 
\proclaim Corollary 6.2. If $f$ is periodic, then $f$ is critical and there is some $n>0$ such that $f_{  Y}^n \Sigma_\gamma = \Sigma_{\beta \gamma}$ and $f_{  Y}^{-n} \Sigma_C = \Sigma_{BC}$.

\noindent{\it Proof.} If $f$ is pseudo-automorphism then so is $f_{  Y}$. Since both $f_{  Y}$ and $f_{  Y}^{-1}$ have the unique exceptional hypersurface, there exists $n \ge 0$ such that $ f_{  Y}^n \Sigma_\gamma = \Sigma_{\beta \gamma}$ which blows up to a hypersurface $\Sigma_C$. If $f$ is periodic then so is $f^{-1}$ and thus $f_{  Y}^{-n} \Sigma_C =\Sigma_{BC}$. \qed

\proclaim Lemma 6.3. Suppose $f$ is periodic. If $f_{  Y}^j \Sigma_\gamma \subset {\cal F}_{0\beta\gamma}$ (possibility (iii) in \S4), then $f_{  Y}^{j+1} \Sigma_\gamma = {\cal F}_{0BC}$.  Thus $f_{  Z}$ is a pseudo-automorphism. 

\noindent{\it Proof.} Suppose $f$ is a periodic map with period $p$. For each $i = 1, \dots, p$, let us set $V_i = f_{  Y}^i \Sigma_\gamma$. It follows that $(f_{  Y}^{-1})^i \Sigma_\gamma = V_{p-i}$. If $V_ j \subset {\cal F}_{0\beta\gamma}$ then $V_{j+1} = {\cal F}_{0BC}$. Applying Lemma 4.10 to $f_{  Y}^{-1}$ we see that $f_{  Y}^{-1}  V_{j+1} = V_{j} = {\cal F}_{0\beta\gamma}$. It follows that $f$ satisfies every condition in Theorem 5.1 and therefore $f_{  Z}$ is a pseudo-automorphism. \qed


\proclaim Lemma 6.4. If $f$ is periodic then $\chi_f(t)$  is self-reciprocal, and $\chi_f=\chi_{f^{-1}}$.  

\noindent{\it Proof.}  A polynomial  $p(z)= \sum_{i=0}^k a_i z^i,\ a_i \in {\bf C}$  is {\it self-reciprocal} if $p(z) =\pm z^k \overline{p(1/\overline{z})}$. If $f$ is periodic then the characteristic polynomial of $f^*_{  Z}$, $\chi(t)$ is a product of cyclotomic factors and thus $\chi(t)$ is self-reciprocal. Furthermore by Lemma 6.3 $f_{  Z}$ is a pseudo-automorphism and therefore $(f_{  Z}^*)^{-1} = (f_{  Z}^{-1})^{*}$. It follows that  $\chi_f$ and $\chi_{f^{-1}}$ are integer polynomials with the same roots. \qed

\proclaim Lemma 6.5. If $f$ is periodic, then there is a non-negative integer $m$ such that 
\itemitem{(1)} $m =m_u = m_d<N$,  $1<d_1 < u_1 < \cdots < d_m < u_m < N$ and 
\itemitem{(2)} $N-u_j = d_{m+1-j}, N-d_j = u_{m+1-j}$ for $j = 1, \dots, m$.

\noindent{\it Proof.} From $(5.8)$ we see that the characteristic polynomial $\chi=\chi_f(t)$ is given by $ \chi (t) = t^{N-1} (t^2+1) \varphi(t)$, where 
$$\eqalign{ \varphi(t) =& {1 \over t^{m_s} }(t-1) (t+1+{1\over t}) + {1 \over t^{N-1}} \left( t^N (t^3-t-1) + (t^3+t^2-1)   \right)  +\cr 
& + (t^3+t^2) \left( t \sum_{j=1}^{m_u} {1 \over t^{u_j}} + (t-1) \sum_{j=1}^{m_d} {1 \over t^{d_j}} \right) - t \sum_{j=1}^{m_u} {1 \over t^{u_j}}-t \sum_{j=1}^{m_d} {1 \over t^{d_j}} .} \eqno{(6.1)}$$
By Lemma 6.4, $\chi(t)$ should be self-reciprocal. Since the first part of $\chi$ and the first line of $(6.1)$ are self-reciprocal, it suffices to consider the case $m_s=  \infty$ and $m_u m_d \ne 0$.  In this case ${\rm dim} f_{  Z}^j \Sigma_\gamma = {\rm dim} f_{  Z}^{j+1} \Sigma_\gamma$ if and only if $j \not\in \{u_i, i=1, \dots, m_u\} \cup \{d_i, i=1, \dots, m_d\}$. Thus it is clear that we have $m =m_u = m_d<N$ and $1<d_1 < u_1 < \cdots < d_m < u_m < N$ for some positive integer $m$. Thus we have
$$ \eqalign{ f_{  Z} : \Sigma_\gamma \to &\Sigma_{BC} \to \cdots \to f_{  Z}^{d_1}\Sigma_\gamma \to p_1 \in \Sigma_{BC} \to \cdots \to q_1 \in \Sigma_{\beta\gamma} \cr & \rightsquigarrow f_{  Z}^{u_1+1} \Sigma_\gamma \subset \Sigma_{C} \to \cdots \to f_{  Z}^{d_2}\Sigma_\gamma \to \cdots \to f_{  Z}^N\Sigma_\gamma = \Sigma_{\beta\gamma}\rightsquigarrow \Sigma_C}$$
By interchanging the roles of $\Sigma_\beta$, $\Sigma_\gamma$ and $\Sigma_B$, $\Sigma_C$, we see that the characteristic polynomial for $f^{-1}$ is given by $ \hat\chi_{f^{-1}}(t) = t^{N-1} (t^2+1) \hat\varphi(t)$ where  
$$\eqalign{ &\hat \varphi(t) ={1 \over t^{N-1}} \left ( t^N (t^3-t-1) + (t^3+t^2-1) \right) + \cr 
& + (t^3+t^2) \left( t \sum_{j=1}^{m} {1 \over t^{N-d_j}} + (t-1) \sum_{j=1}^{m} {1 \over t^{N-u_j}} \right) - t \sum_{j=1}^{m} {1 \over t^{N-d_j}}-t \sum_{j=1}^{m} {1 \over t^{N-u_j}} .} \eqno{(6.2)}$$ 
Since both $f$ and $f^{-1}$ have the same characteristic polynomial, we obtain the second statement of the Lemma by comparing  $\chi_f$ and $\chi_{f^{-1}}$.  \qed



\proclaim Lemma 6.6. Suppose $f$ is periodic. 
\itemitem{(a)} If $m$ is even, then for all $j = 1, \dots, m$, $2 \le u_{j}-d_j \le d_1$. 
\itemitem{(b)} If $m$ is odd, then $1 \le u_{(m+1)/2}-d_{(m+1)/2} \le d_1$ and for all $j \ne (m+1)/2$, $2 \le u_{j}-d_j \le d_1$.

\noindent{\it Proof.} Suppose $j_*$ is the smallest positive integer such that $u_{j_*} - d_{j_*} > d_1$. Then we have (1) $f_Y^{d_{j_*}} \Sigma_\gamma = $ a fibration in $\Sigma_\gamma$, (2) $f_Y^{d_{j_*}+i} \Sigma_\gamma \in {\rm\ a \ point\ in\ } \ f_Y^{i} \Sigma_\gamma$ for $i = 1, \dots, d_1$, and (3) $f_Y^{d_{j_*} + d_1+1} \Sigma_\gamma \in f_Y^{d_1+1} \Sigma_\gamma$ which is a point in $\Sigma_{BC}$. It follows that the exceptional hypersurface $\Sigma_\gamma$ is pre-periodic which contradicts to the hypothesis $f$ is periodic. If  $ u_{j}-d_j=1$ then $f^{d_j+1}_Y \Sigma_\gamma = \Sigma_{BC} \cap \Sigma_{\beta\gamma} = f^{u_j}_Y \Sigma_\gamma$. Thus the situation  $ u_{j}-d_j=1$ can happens at most once and by Lemma 6.5 we see that $(N-d_j) -(N-u_{j}) = u_{m-j+1} - d_{m-j+1} =1$. It follows that $j = m-j+1$ and thus $j = (m+1)/2$. \qed

\proclaim Lemma 6.7. Suppose $f$ is critical and $m \ge 1$ then
\itemitem{(a)} $d_1\ne 1,3, 4$
\itemitem{(b)} If $d_1 =2$ then $m=1$, and  either (i) $\alpha_0=\alpha_2=1$ and $\beta_0 =0$ or (ii) $\alpha_0= \eta^2 \alpha_2 = \eta, \beta=\eta^2$ where $\eta^2 -\eta+1=0$.
\itemitem{(c)} If $m\ge 2$ is odd then for $j= 1, \dots, m-1$, $d_{j+1} - u_j \ge 5$. 
\itemitem{(d)} If $m\ge 2$ is even then for $ 1\le j\le m-1$ and $j \ne m/2$, $d_{j+1} - u_j \ge 5$ and $d_{m/2+1} - u_{m/2} \ge 4$. 

\noindent{\it Proof.} (a) $d_1=1$ means $\Sigma_{BC}$ is a line through $e_1$ in $\Sigma_\gamma$. Since $\Sigma_{BC} = \{ x_3=0, \alpha_0 x_0+ \alpha_2 x_1 + x_2=0\}$ and $\alpha_2\ne 0$, it follows that $e_1 \not\in \Sigma_{BC}$ and thus $d_1 \ne 1$. Since $\Sigma_{BC} \subset \Sigma_3$, we have $f_Y^3 \Sigma_\gamma = f_Y^2 \Sigma_{BC} \subset \Sigma_1$ which doesn't contain $e_1$. Furthermore $f_Y\Sigma_1 = \{ [\beta_0 x_0: \beta_0x_2: \beta_0 x_3: \alpha_0 x_0 + \alpha_2 x_2 + x_3\}$ if $\beta_0 \ne 0$  and  $f_Y\Sigma_1$ is a line in the blowup fiber $E_3$ if $\beta_0 = 0$. It follows that $f_Y^4 \Sigma_\gamma$ does not contain $e_1$. Therefore $d_1 \ne 3$ or $4$. The statement for (b) can be confirmed by direct computation. For each $j \le m-1$, $f_Y^{u_j+1} \Sigma_\gamma$ is a line through $e_3$ in $\Sigma_C$ which can be parametrized as $t\mapsto \{ [ 1: \mu: -\alpha_0 - \alpha_2 \mu: t]\}$ for some fixed $\mu\in {\bf C} \cup \{ \infty\}$. By computing the forward iteration of $[ 1: \mu: -\alpha_0 - \alpha_2 \mu: t]$ we see that $d_{j+1} - u_j \ne 1,2,$ or $3$. Furthermore $d_{j+1} - u_j = 4 $ if and only if $f_Y^{u_j+1} \Sigma_\gamma = \Sigma_C \cap \{ (1+ \alpha_2 \beta_0) x_0 + \alpha_2 x_1 =0 \}$. It follows that $d_{j+1}-u_j=4$ occurs only once. Suppose  $d_{j+1}-u_j=4$ for some $1 \le j\le m-1$. By Lemma 6.5 we see that $(N-u_j) -(N-d_{j+1}) = d_{m-j+1} - u_{m-j} =4$. It follows that $j+1 = m-j+1$ and thus $j = m/2$. The statement (c) and (d) follows. \qed

Direct computation shows the following properties: 
\proclaim Lemma 6.8. Suppose $f$ is critical, then $\hat\varphi(t)$ defined in $(6.2)$ satisfies
 \itemitem{(a)} $\hat\varphi(1) =0$, and 
 \itemitem{(b)} $\hat\varphi'(1) = 7(m+1) - (N +\sum_{j=1} ^m (u_{j}-d_j))$.

\proclaim Lemma 6.9. Suppose that $m\ge 2$ and that $f$ is critical satisfying (5.1).  Then $$ N + \sum_{j=1}^{m} (u_j - d_j ) \quad\ \   \eqalign{&\ge\ \ 9m + 3 \qquad {\rm \ \ if\ \ }m{\rm \ is\ odd} \cr& \ge \ \ 9m + 4 \qquad {\rm \ \ if\ \ }m{\rm \ is\ even} \cr}$$
Thus $\hat\varphi'(1)<0$, so $\hat\varphi$ has a root greater than 1.

\noindent{\it Proof.} Suppose $m$ is even. By Lemma 6.5 we see that
$$\eqalign{N + \sum_{j=1}^{m} (u_j - d_j ) = &2 d_1 + 4(u_1-d_1) + 2 ( d_2 - u_1) + \cdots\cr &\cdots+ 2( d_{m/2} - u_{m/2-1}) + 4(u_{m/2} - d_{m/2}) + (d_{m/2+1}-u_{m/2})\cr}$$ 
By Lemma 6.7, (a) and (b), we have $d_1\ge5$.  Applying Lemma 6.6 (a) and Lemma 6.7 (d) we have $$N + \sum_{j=1}^{m} (u_j - d_j )  \ge 2 \cdot 5 + 4\cdot 2 + 2\cdot 5 + \cdots + 4 \cdot 2 +  4 = 9m+4.$$
Similarly when $m$ is odd 
$$\eqalign{N + \sum_{j=1}^{m} (u_j - d_j ) = &2 d_1 + 4(u_1-d_1) + 2 ( d_2 - u_1) + \cdots\cr &\cdots+ 2( d_{(m+1)/2} - u_{(m-1)/2}) + 2(u_{(m+1)/2}-d_{(m+1)/2}).\cr}$$
Again applying Lemma 6.6 (b) and Lemma 6.7(c) we have 
$$N + \sum_{j=1}^{m} (u_j - d_j )  \ge 2 \cdot 5 + 4\cdot 2 + 2\cdot 5 + \cdots +2 \cdot 5 +  2 = 9m+3.$$
\qed

\proclaim Theorem 6.10.   If $f$ is periodic with  $m=0$ and $m_s=\infty$ then $f$ is one of the following: 
\item{$\bullet$ } $\alpha=(-1,0,-1,1),\ \beta=(0,1,0,0)$ :  $f_{\alpha\beta}$ has period $8$ and there is a conic $Q$ such that $$f_{  Y}: \Sigma_\gamma \to \Sigma_{BC} \to Q \to \Sigma_{\beta\gamma} \rightsquigarrow \Sigma_C.$$
\item{$\bullet$ } $\alpha=(-1/2,0,-1,1), \ \beta=(1,1,0,0)$:  $f_{\alpha\beta}$ has period $12$, and $$f_{  Y}: \Sigma_\gamma \to \Sigma_{BC} \to L_1 \to L_2 \to \Sigma_{\beta\gamma} \rightsquigarrow \Sigma_C$$
where we set $L_1 = \Sigma_2 \cap \{x_0+x_3=0\}$ and $L_2= \Sigma_1 \cap \{x_0+x_2=0\}$. 

\noindent{\it Proof.}  The polynomial defined in (5.8) is also given by 
$$ \chi(t) = (t^2+1) \left(t^N (t^3-t-1) + t^3+t^2 -1 \right).$$ 
It follows that $\chi(t)$ has a root bigger than $1$ if and only if $N\ge 8$ and in case $N=7$ the matrix representation of $f_{  Z}^*$ has $3\times 3$ Jordan block with eigenvalue $1$.  Thus we need to check the situation $f^{n+1}\Sigma_\gamma=\Sigma_{\beta\gamma }$ only for $n\le 5$.  For this, let us parametrize $\Sigma_{BC} = \{ [ 1: t: -\alpha_0 - \alpha_2 t:0]\}$ and let $[f_0^{(n)}:f_1^{(n)}:f_2^{(n)}:f_3^{(n)}]$ denote the $n$-th iteration of $\Sigma_{BC}$. If $ f^{n+1} \Sigma_\gamma = \Sigma_{\beta\gamma}$ then for all $t$ we have 
$$\beta_0 f_0^{(n)}+f_1^{(n)} =0,\ \ \ {\rm and} \ \ \ \alpha_0 f_0^{(n)}+ \alpha_2 f_2^{(n)}+ f_3^{(n)}=0 \eqno{(6.3)}$$
Since equations in (6.3) are polynomials in $t$ whose coefficients are integer polynomials in the variables $\beta_0, \alpha_0,$ and $\alpha_2$,  we may use the  computer show that for $0\le n\le 5$, the only two possibilities are those listed above.  \qed

\proclaim Theorem 6.11. If $f$ is periodic with  $m=1$ and $m_s=\infty$,  then $f$ is one of the following:
\item{$\bullet$ } $\alpha=(1,0,1,1),\ \beta=(0,1,0,0)$ :   $f$ has period $8$,  $\Sigma_{BC} \cap \Sigma_{\beta\gamma} \ne \emptyset$, and 
$$f_{  Y}: \Sigma_\gamma \to \Sigma_{BC} \to \Sigma_\gamma \cap \Sigma_2 \to \Sigma_{BC} \cap \Sigma_{\beta\gamma}  \rightsquigarrow \Sigma_C\cap \Sigma_2 \to \Sigma_{\beta\gamma} \rightsquigarrow \Sigma_C.$$
\item{$\bullet$ } $\alpha=(\eta /(1-\eta),0,\eta,1), \ \beta=(\eta^2,1,0,0)$ and $\eta^3=-1$, $\eta \ne -1$:  $f_{\alpha\beta}$ has period $12$, and
$$f_{  Y}: \Sigma_\gamma \to \Sigma_{BC} \to \Sigma_\gamma \cap \Sigma_2 \to p_1 \in \Sigma_{BC} \to p_2 \in \Sigma_{\beta\gamma}  \rightsquigarrow \Sigma_C\cap \Sigma_1 \to \Sigma_{\beta\gamma} \rightsquigarrow \Sigma_C,$$
where $p_1 = [1:0:-\eta^2:0] \in \Sigma_{BC}$ and $p_2 = [1:-\eta^2:0:-\eta^2] \in \Sigma_{\beta\gamma}$.

\noindent{\it Proof.} From (5.8) the characteristic polynomial of $f^*_{  Z}$ is given by 
$$ \chi(t) =t^{N-(u_1+d_1)} (t^2+1) (t^{d_1}+1) \left (t^{u_1} (t^3-t-1) + t^3+t^2 -1\right ).$$ 
It follows that $\chi(t)$ has a root bigger than $1$ if and only if $u_1 \ge 8$. If $u_1=7$, the $f_{  Z}^*$ has a $3\times 3$ Jordan block. Thus If $f_{  Z}^*$ is  periodic then $d_1 \le5 < u_1$. By direct computation of $f^n\Sigma_\gamma = f^{n-1} \Sigma_{BC}$ for $n=1, \dots, 5$, we can easily check the two conditions (i) $f^{n-1} \Sigma_{BC} \subset \Sigma_\gamma$, (ii) $f^{n-1} \Sigma_{BC} \subset \{ x_3 = \lambda x_2\} $ for some $\lambda \in {\bf C}$ and thus we see that there are only two possibilities listed in this Theorem.\qed

\proclaim Theorem 6.12. If $m\ge 2, m_s=\infty$ then $f$ has exponential degree growth (and is not periodic).

\noindent{\it Proof.} By Lemmas 6.8 and 6.9 we see that $\chi_N(1) = 0$ and $\chi_N'(1) = 2 \hat \varphi'(1) <0$. Since the leading coefficient of $\chi_N$ is $1$, there exist a real root which is  strictly bigger than $1$. It follows that the dynamical degree of $f$ is strictly bigger than $1$.  \qed

\proclaim Theorem 6.13.  If $1\le m_s <\infty$, then $f$ is not periodic.


\noindent{\it Proof.}  By Lemmas 6.8 and 6.9 we see that $\chi_N(1) = 0$ and $\chi_N'(1) = 2(3+ \hat \varphi'(1) )=2(-2 m+1)$. It follows that if $m\ge 1 $ then $f$ has positive entropy. Now suppose $m=0$, we have 
$\chi_N(t) =  (t^2+1) \left(t^N(t^3-t-1) + t^{N-m_s-2}(t-1) (t^2+t+1) + t^3+t^2-1\right).$ 
If $f$ is periodic then the characteristic polynomial for $f_{  Z}^*$ should be self-reciprocal. It follows that $m_s = (N-4)/2$. It follows that $N$ should be even and $$\chi_N(t) =   (t^2+1) \left(t^N(t^3-t-1) + t^{N/2}(t-1) (t^2+t+1) + t^3+t^2-1\right).$$ 
By inspection we see that $m_s \ge 3$ and it follows that $N\ge 10$. We can also check that $\chi_N(1) =0$ and $\chi'_N(t) = 10-N$. Therefore if $m_s \ne 3$ then $f$ is not periodic. In case $m_s=3$, the matrix representation for $f_{  Z}^*$ has $3\times3$ Jordan block with eigenvalue $1$ and all other eigenvalues have modulus $1$. \qed 

\medskip\noindent{\it Proof of  Theorem 5: }  The statement of the Theorem 5 in the Introduction follows from Theorems 6.10--13. \qed

We remark that in the proof of Theorem 6.13, we see that if $m_s=3$ and $m=0$, then the degree of $f^n$ is quadratic in $n$.  This case occurs for   $\alpha=(a,0,1,1)$ and $\beta=(0,1,0,0)$, which is the so-called Lyness process and will be discussed in \S8.

\bigskip

\noindent{\bf\S 7. Pseudo-automorphisms with positive entropy} In this section we consider the case
$$\beta=(0,1,0,0) \quad {\rm and} \quad \alpha=(a,0,\omega,1) \eqno{(7.1)}$$ where $\omega^2+\omega+1=0$ and $a \in {\bf C}\setminus\{0\}$. With this choice of parameters, we see that $f$ is critical and that $\Sigma_B=\Sigma_3$ and $\Sigma_\beta=\Sigma_1$.  Since the maps $f:\Sigma_3\to\Sigma_2\to \Sigma_1$ are dominant, (4.3) gives an 8-cycle of dominant maps
$$ f_{  Y} : \Sigma_1 \to E_3\to S_{01} \to \Sigma_0 \to S_{03} \to E_1 \to \Sigma_3 \to \Sigma_2 \to \Sigma_1 \eqno(7.2)$$
Since this 8-cycle is fundamental to our understanding of $f$ in this case, we will refer to the union of these 8 hypersurfaces as the {\it rotor} and denote it as ${\cal R}$.  Clearly, $f^8_{  Y}$ fixes each component of the rotor; in addition, it  has a relatively simple expression.  On $\Sigma_3$, or example, we have:
$$\eqalign{ f^8_{  Y} : \Sigma_3 \ni & [x_0: x_1: x_2:0] \cr & \mapsto [ x_0 ( a x_0 + \omega x_2) (a x_0+a x_1 + \omega x_2) \cr
&:x_1 ( x_1 x_2 + a \omega x_0^2 + a \omega x_0 x_1 + a \omega x_0 x_2 + \omega^2 x_0  x_2 + \omega^2 x_2^2)\cr
&:\omega x_2 (a x_0 + \omega x_2) (x_1+ a \omega x_0 + \omega^2 x_2) : 0] \in \Sigma_3. \cr }\eqno(7.3)$$
The restriction of $f_Y^8$ to the rotor is studied in \S{\bf C}.

Note that by (7.1), $\Sigma_{BC}=\Sigma_3\cap \Sigma_C$ and $\Sigma_{\beta\gamma}=\Sigma_1\cap\Sigma_{\gamma}$.  Using (7.2) we may verify that $f_{  Y}$ satisfies condition (5.1), which in this case is
$$f^j_{  Y}\Sigma_\gamma\not\subset {\cal F}_{0\beta\gamma} {\rm\ for\ all\ }1\le j\le 10, {\rm\ and\ } f^{11}_{  Y}=\Sigma_{\beta\gamma} \eqno(7.4)$$
We define the space $\pi_{  Z} : {  Z} \to {  Y}$ by successively blowing up the $11$ curves $\gamma_{j}:=f^j_{  Y} \Sigma_\gamma$, $1\le j\le  11$.  The dynamical degree, being a birational invariant, is independent of the order in which the  $\gamma_j$'s are blown up.

\proclaim Theorem 7.1. The induced map $f_{  Z}$ is a pseudo-automorphism, and the dynamical degree of $f$ is greater than 1. 

\noindent {\it Proof.} From (7.4) we see that $f_{Y}$ satisfies the condition in Theorem 5.1, so $f_Z$ is a pseudo-automorphism.  By Theorem 5.3, the characteristic polynomial of $f^*_Z$ is $t^{11}(t^3-t-1)+t^3+t^2-1=(-1 + t) (1 + t) (1 + t^4) (1 - t^3 - t^4 - t^5 + t^8)$.   Thus $\delta(f)$ is the largest root of this polynomial, which is approximately 1.28064.   \qed

The space $Z$ has been defined earlier, but now let us be more precise:  we define $Z$ as the space obtained by blowing up first $\gamma_{11}\subset Y$, then we blow up the strict transform of $\gamma_{10}$ in the resulting space, followed by blowing up the strict transform of $\gamma_9$, and continuing this way until we blow up the strict transform of $\gamma_1$.      We will use the notation $\Gamma_j$ to denote the exceptional divisor of the blowup of $\gamma_j$.  There are no points where three distinct $\gamma_j$'s intersect.  If $p=\gamma_j\cap\gamma_k$, with $j > k$, then we blow up $\gamma_j$ first, and we refer to the fiber in $\Gamma_j$ over $p$ as the {\it first} fiber over $p$, and write it as ${\cal F}^1_p$.  We then blow up the strict transform of $\gamma_k$, and the blowup fiber over the point $\gamma_k\cap\Gamma_j$ is equal to $\Gamma_j\cap\Gamma_k$.   
  
Let us describe some of the intersections of the $\gamma_j$'s.  $f$ is constant on each line in $\Sigma_\gamma$ passing through $e_1$.  Further, $\gamma_1\subset\Sigma_3$, $\gamma_2\subset\Sigma_2$, and $\gamma_6\subset\Sigma_0$, and $e_1=\Sigma_0\cap\Sigma_2\cap\Sigma_3$.  We set $\ell_2=\Sigma_\gamma\cap\Sigma_3$, $\ell_3=\Sigma_\gamma\cap\Sigma_2$, and $\ell_7=\Sigma_\gamma\cap\Sigma_0$.    Thus we have $f(\ell_j)=\gamma_1\cap\gamma_j$ for $j=2,3,7$.  The curve $\gamma_9\subset\Sigma_1$ is a conic, and $\gamma_9\cap\gamma_1$ consists of two points.  We let $\ell_9'$ and $\ell_9''$ denote the two lines in $\Sigma_\gamma$ for which $f(\ell_9'\cup\ell_9'')=\gamma_1\cap\gamma_9$.  This accounts for all the curves $\gamma_j$ which intersect $\gamma_1$.     As a consequence of the order of blowup, the first fiber ${\cal F}^1_{f(\ell_j)} = {\cal F}^1_{\gamma_{1}\cap\gamma_j}$,  $j=2,3,7$ is contained in $\Gamma_j$ and similarly for $\ell'_9, \ell''_9$. 
  
There is a similar situation for the $\gamma_j$'s which intersect $\gamma_{11}$.    The curves  $\gamma_5$, $\gamma_9$ and  $\gamma_{10}$ each intersect $\gamma_{11}$ in a single point, and $\gamma_3$, which intersects $\gamma_{11}$ in 2 points, and this accounts for all the intersection points between $\gamma_{11}$ and the other $\gamma_j$'s.

Let us use the notation $\pi_1:Z_1\to Y$  for the manifold obtained by blowing up the curve  $\gamma_{11}\subset Y$.  This is the first blowup performed in the construction of $Z$.  Let $f_{Z_1}:Z_1\dasharrow Z_1$ be the induced map.   Since ${\cal I}(f_Y)=\gamma_{11}\cup\Sigma_{02}\cup {\cal F}^1_{0\beta\gamma}$, it follows that ${\cal I}(f_{Z_1})\subset \Sigma_{02}\cup{\cal F}^1_{0\beta\gamma}\cup\gamma_{10}\cup\Gamma_{11}$.

%
%
%
%

\proclaim Lemma 7.2.  ${\cal I}(f_{Z_1})={\cal F}_{0\beta\gamma}^1\cup\Sigma_{02}\cup \gamma_{10}$.  

\noindent {\it Proof. }  We have seen already that the indeterminacy locus is contained in $ \Sigma_{02}\cup{\cal F}^1_{0\beta\gamma}\cup\gamma_{10}\cup\Gamma_{11}$, so it suffices to show that ${\cal I}(f_{Z_1})\cap\Gamma_{11}$ consists of the two points $\gamma_{10}\cap \Gamma_{11}$ and ${\cal F}^1_{0\beta\gamma}\cap\Gamma_{11}$.  Thus we look at $f_{Z_1}$ in coordinate charts that cover $\Gamma_{11}$.  We will look first at $\Gamma_{11}\cap\pi^{-1}(\gamma_{11}-\gamma_{10})$.

In the local coordinates $(s, \zeta, x_3)_{S_{01}} \mapsto [s: s \zeta: 1: x_3] \in {\bf P}^3$ in the neighborhood of $S_{01} - E_1= \{ s=0, \zeta\ne \infty\}$, we have $\gamma_{11} = \{ \zeta =0, a s + \omega + x_3=0\}$ and ${\cal F}^1_{0\beta\gamma} \cap \gamma_{11} = (0,0,-\omega)_{S_{01}}$. 
We use the local coordinate charts $(s,t,\eta)'$  on $U'$ and $(s,\eta,t)''$ on  $U''$ so that $\pi_1$ is given by 
$$\eqalign{ &\pi':  U'\ni (s, t, \eta)' \mapsto (s,t, -a s-\omega+t \eta)_{S_{01}}\cr
&\pi'': U''\ni (s, \eta,t)'' \mapsto (s,t\eta, -a s-\omega+t )_{S_{01}}\cr} $$ 
It is evident that $\Gamma_{11}\supset\{t=0\}$ in both coordinate charts, and $U'\cup U'' \supset \pi_1^{-1} (\gamma_{11} - \gamma_{10})$.
The induced map $f\circ\pi_1:U'\cup U''\to{\bf P}^3$ is given by 
$$ \eqalign{ &\ \  U'\ni (s,t, \eta)' \mapsto [s:1:-as+t \eta-\omega:\eta]\ \  \cr
&\ \  U''\ni (s,\eta, t)'' \mapsto [s \eta:\eta:\eta(-as+t-\omega):1]  .\cr}\eqno{(7.4)}$$
so we see that $\{t=0\}$ is mapped to $\Sigma_C$.

From $(7.4)$ we see that the map $f\circ\pi_1:U'\cup U''\to {\bf P}^3$ is everywhere regular.  The only points of $\Sigma_C$ which is blown up in the construction of $Y$ are $e_3$ and $[0:1:-\omega:0]$  which is point the base point of ${\cal F}^1_{0BC}$.   By (7.4), the preimage of $e_3$ is $(0,0,0)''\in U''$, and the preimage $[0:1:-\omega:0]$ is $(0,0,0)'\in U'$.  Working in local coordinates in $Y$ over $e_3$, we find that $f\circ\pi_1:U''\to Y$ is everywhere regular.    Thus we conclude that $ f_Y\circ \pi_1:Z_1\to Y$ is regular on $(U'-(0,0,0)') \cup U''$.   Now in order to pass to $f_{Z_1}$ we need to consider the point $\gamma_{11}\cap\Sigma_C$ which is blown up.  However, this is the image point of $\gamma_{10}\cap\gamma_{11}$, which is not in our coordinate chart.   We note that $(0,0,0)'$ is the point ${\cal F}^1_{0\beta\gamma} \cap \Gamma_{11}$, so we conclude that $f_{Z_1}$ is regular at all points of $\Gamma_{11}-(\pi_1^{-1}(\gamma_{11}\cap\gamma_{10})\cup{\cal F}^1_{0\beta\gamma})$.


Now we consider $\Gamma_{11}\cap \pi^{-1}(\gamma_{11}-{\cal F}^1_{0\beta\gamma})$, which does not lie over any of the centers of blowup in the construction of $Y$. As in (5.7), we may use the local coordinates  $(s, x_2, \zeta) \mapsto [1:s:x_2:-a-\omega x_2 + s \zeta]$ in a neighborhood of $\{s=0,\zeta\ne\infty\}\subset \Gamma_{11}-\pi_1^{-1}({\cal F}^1_{0\beta\gamma})$, and we get
$$f_{Z_1} : \Gamma_{11} \ni  (0,x_2,\zeta) \mapsto [1:x_2:-a-\omega x_2:\zeta]\in \Sigma_C \qquad {\rm if\ } (0,x_2,\zeta)  \ne (0,0,a \omega-a).\eqno{(7.5)}$$
Similarly using the local coordinates $(\zeta, x_2, s) \mapsto [1:s \zeta :x_2:-a-\omega x_2 + s] \in \Gamma_{11}$, we have  
$$f_{Z_1} : \Gamma_{11} \ni  (\zeta,x_2,0) \mapsto [\zeta:x_2\zeta:\zeta(-a-\omega x_2): 1]\in \Sigma_C \ \  {\rm if\ } (\zeta,x_2,0)  \ne ({1\over a \omega-a},0,0).\eqno{(7.6 )}$$

Since both $(0,0, a\omega-a)$ in $(7.5)$ and $(1/( a \omega-a),0,0)$ in  $(7.6)$ correspond to the point $\gamma_{11} \cap \gamma_{10}$, combining with the previous conversation about $\Gamma_{11} - \pi_1^{-1}(\gamma_{11}\cap\gamma_{10})$, we conclude that $f_{Z_1}$ is regular at all points of $\Gamma_{11} - (\gamma_{10} \cup {\cal F}^1_{0\beta\gamma})$.
\qed

\proclaim Lemma 7.3.  The three curves $\gamma_5$, $\gamma_{11}$, ${\cal F}^1_{0\beta\gamma}$  intersect transversally inside $Y$, and $\gamma_1$, $\gamma_7$, ${\cal F}^1_{0BC}$ intersect transversally inside $Y$.  Thus, inside $Z_1$, the strict transform of ${\cal F}^1_{0\beta\gamma}$ is disjoint from the strict transforms of $\gamma_j$, $1\le j\le 10$. 

\noindent{\it Proof.}  It suffices to prove the first statement.  We may write $\gamma_{11}\subset{\bf P}^3$ as $s\mapsto [s:0:1:-as-\omega]$.  This intersects $\Sigma_{01}$ in the point $[0:0:1:-\omega]$.  We use the coordinate system $\pi:(u,\eta,x_3)\mapsto [u:u\eta:1:x_3]\in{\bf P}^3$.  Thus ${\cal F}^1_{0\beta\gamma}=\{u=0, x_3=-\omega\}$.  In this coordinate system,  $\gamma_{11}$ becomes $s\mapsto (s,0,-as-\omega)$, so $\gamma_{11}$ crosses $S_{01}$ when $s=0$, at the point $(0,0,-\omega)$.    On the other hand, if we map $\gamma_{11}$ backward under $f_Y^{-6}$, we find an expression for $\gamma_5$.  The base point is given by $[0:0:1:s]$, and the fiber coordinate is given by $\eta=(1+as)(1+as+\omega(1+a - s))/(as (-1+s +as -\omega(1+a-as))$.  Thus when the base point is $[0:0:1:-\omega]$, we have $\eta=0$.  Thus all three curves meet at $(u,\eta,x_3)=(0,0,-\omega)$.  The curve $\gamma_{11}$ is transverse to $S_{01}$,  but $\gamma_5$ and ${\cal F}^1_{0\beta\gamma}$ are tangential to $S_{01}$, so $\gamma_{11}$ is transverse to the other two, and $\gamma_5$ is transverse to $\{x_3=-\omega\}$, while ${\cal F}^1_{0\beta\gamma}$ is tangential to this set.  \qed

 For $2 \le j\le 11$, let $\pi_j : Z_j \to Z_{j-1}$ be the blowup of the strict transform of $\gamma_{12-j}$ inside $Z_{j-1}$ and $\pi: Z \to Y = \pi_{11} \circ \pi_{10}\circ \cdots \circ \pi_1$, that is, we blowup $\gamma_{11}$ first, then $\gamma_{10}$, then $\gamma_{9}$, etc.  Let $f_{Z_j}:Z_j \dasharrow Z_j$, $f_Z: Z \dasharrow Z$ denote the induced map. 

\proclaim Lemma 7.4.  For $1\le j \le 10$, ${\cal I}(f_{Z_j}) = {\cal F}^1_{0\beta\gamma} \cup \Sigma_{02} \cup \gamma_{11-j}$.

\noindent {\it Proof. }  Suppose $p$ is a point of $\gamma_j \cap \gamma_k$ $1 \le j<k\le 10$. Because of the order of blowup, $\gamma _k$ is blown up before $\gamma_j$ and $\gamma_{k+1}$ is blown up before $\gamma_{j+1}$. Since $f_Y$ is regular at $p$ and the order of blowups at $p$ is consistent with the order of blowups at $f_Y(p)$,  the induced map $f_{Z_i}$ is a local biholomorphism in a neighborhood of the exceptional divisor over $p$ for $12-j\le i \le 11$.

Notice that for all $1\le j\le 11$ the strict transformation of $\gamma_j$ does not intersect $\Sigma_{02}$ in $Y$. Suppose $\gamma_j$ intersects ${\cal F}^1_{0\beta\gamma}$ at a point $q$. Using the local coordinates in the neighborhood $(s, \zeta, x_3)_{S_{01}}$, we may assume that $q= (0, \zeta_*, - \omega)_{S_{01}}$ and $\gamma_j(s) = (Q_1(s), Q_2(s)+ \zeta_*, Q_3(s)-\omega)_{S_{01}}$, where $\gamma_j = \{ \gamma_j(s), s \in {\bf C}\},$ and $ \gamma_j(0) = q$. Consider two local coordinate charts covering the exceptional divisor over the point $q$ : 
$$\eqalign{ &( s, t, \eta) \mapsto ( Q_1(s), Q_2(s) + \zeta_* + t, Q_3(s) - \omega + t \eta)_{S_{01}} \cr &( s, \eta, t) \mapsto ( Q_1(s), Q_2(s) + \zeta_* + t \eta, Q_3(s) - \omega + t )_{S_{01}} \cr }$$
With a computation similar to Lemma 7.2, We see that the induced map is regular everywhere on the exceptional divisor over q, ${\cal F}(q)$, except the point of intersection ${\cal F}(q) \cap {\cal F}^1_{0\beta\gamma}$. Now since the curve $\gamma_{11-j}$ is the pre-image of $\gamma_{12-j}$, we have ${\cal I}(f_{Z_j}) = {\cal F}^1_{0\beta\gamma} \cup \Sigma_{02} \cup \gamma_{11-j}$. \qed

From the previous Lemma we have ${\cal I}(f_{Z_{10}}) = {\cal F}^1_{0\beta\gamma} \cup \Sigma_{02} \cup \gamma_{1}$. Since $\Sigma_\gamma$ is the pre-image of $\gamma_1$, we have ${\cal I}(f_Z) \subset  {\cal F}^1_{0\beta\gamma} \cup \Sigma_{02} \cup \Sigma_\gamma$. From $(5.2)$ we see that for all most every line $\ell\subset\Sigma_\gamma$, through $e_1$ in $\Sigma_\gamma$, $f$ maps $\ell$ regularly to a point $q\in\gamma_1$.   In our construction of $Z$, we blew up $\gamma_{11}, \dots, \gamma_2$ before $\gamma_1$.  Thus the map $f_Z$ will map $\ell$ regularly to the fiber of $\Gamma_1$ over $q$ unless $q$ is an intersection point of $\gamma_1\cap \gamma_j$ for some $2\le j\le 11$.

\proclaim Lemma 7.5. Suppose $q \in \gamma_1 \cap \gamma_j$ for some $j=2, \dots, 11$ and $\ell_j \subset \Sigma_\gamma$ be the line which mapped to $q$ by $f_Y$. The line $\ell_j \subset {\cal I}(f_Z)$ and every point in $\ell_j$ blows up to the first blowup fiber ${\cal F}^1_q$.  

\noindent{\it Proof.} Let us parametrize $\gamma_1 = \{ \gamma_1(t) = [ -{1 \over a} (1+ \omega t): t: 1:0], t \in {\bf C} \}$. Let us set $q = \gamma_1( t_*)$ for some $t_* \in {\bf C}$ and $\gamma_j = \{ \gamma_j(s) = [ Q_0(s)  -{1 \over a} (1+ \omega t_*):Q_1(s)+ t_*: Q_2(s)+1:Q_3(s)]\}$. The line $\ell_j$ is given by the strict transform in $Y$ of the line connecting $e_1$ and $\tilde q =  [ -{1 \over a} (1+ \omega t_*):0: t_*: 1]$ in ${\bf P}^3$. To see the image of the line $\ell_j$, we  consider the set $U = \{  [ -{1 \over a} (1+ \omega t_*) + s \zeta :u : t_*+s: 1]\}$ which has the property that $U \cap \{ s=0\}  = \ell_j - \{ e_1\}$. 
Since the point $q$ is blown up twice, let us consider a local coordinate charts for $\pi_{12-j}^{-1} ( \gamma_j )$ : $$( v, \xi, s)_{\gamma_j} \mapsto  \left [{ Q_0(s)  -{1 \over a} (1+ \omega t_*)\over Q_2(s) +1} + v:{Q_1(s)+ t_*\over Q_2(s) +1}+ v \xi: 1:{Q_3(s)\over Q_2(s) +1}\right].$$
Using the induced map $f_Z$, we see that 
$$f_Z: \ell_j \ni [ -{1 \over a} (1+ \omega t_*)  :u : t_*: 1] \rightsquigarrow \{ (0,\xi,0), \xi \in {\bf C}\} \subset \Gamma_j, $$
that is, each point in $\ell_j$ blows up to a whole first blowup fiber over $q$. \qed

Before Lemma 7.2, we enumerated the possibilities for lines $\ell$ and points $q$ as in the hypotheses of Lemma 7.5.  Thus we may combine Lemmas 7.2--5 to have the following Theorem: 

\proclaim Theorem 7.6. The indeterminacy locus ${\cal I}(f_Z) = \Sigma_{02} \cup {\cal F}^1_{0\beta\gamma} \cup \ell_2 \cup \ell_3 \cup \ell_7 \cup \ell'_9 \cup \ell''_9$.   If $\zeta$ is a point of one of the lines $\ell$, then $f_Z$ blows up $\zeta$ to the first fiber ${\cal F}^1_{f(\ell)}$.

Now we give the existence of Green currents for the invariant class $\alpha=\alpha_Z^+\in H^{1,1}(Z)$.  

\proclaim Theorem 7.7.  There is a positive closed current $T^+_Z$ in the class of $\alpha^+_Z$ with the property: if \ $\Xi^+$\/ is a smooth form which represents $\alpha^+_Z$, then $\lim_{n\to\infty}\delta_1(f)^{-n} f_{Z}^{ n*}\Xi_Z^+ =T^+_Z$  in the weak sense of currents on $Z$.

\noindent{\it Proof.}  The map $f_Z^*$ is given in Appendix A, where we are in Case (II).  Working directly with the matrix (A.1), we see that the invariant class is given by:
$$\alpha=H_Z - c_1\tilde E_1-c_3 \tilde E_3-c_{01}\tilde S_{01}-c_{03}\tilde S_{03} - \sum_{j=1}^{11}c'_j {\cal F}_j $$
where $c_1,c_3>0$, $c_1+c_3=1$, $c'_{11}>c'_{10}>\cdots>c'_1>0$, and  $c_{01}=c_{03}>c'_8$.  As in Theorem~4.9, we will show that $\alpha^+_Z\cdot\sigma$ for each curve $\sigma$ inside the forward image of ${\cal I}(f_Z)$.  The result will the follow from Theorem 1.3 of [Ba].

   Let us start with ${\cal F}_{0\beta\gamma}\subset{\cal I}(f_Z)$.  Points of this curve are blown up to ${\cal F}_{0BC}$.  The curve $\sigma={\cal F}_{0BC}$ is the exceptional fiber inside $S_{03}$ over the point $\Sigma_{BC}\cap \Sigma_{03}\in {\bf P}^3$.   Thus  $\sigma\cdot S_{03}=-1$.     In the construction of $Z$,  $\gamma_7$ will be blown up to create the exceptional divisor $\Gamma_7$.  At this stage, by Lemma 7.3, $\sigma$ and $\gamma_1$ become separated.   Thus $\sigma\cdot\Gamma_1=0$, and $\sigma\cdot\Gamma_7=1$, so , so $\alpha\cdot\sigma=c_{03}-c_7>0$.  
   
Points of the indeterminate curve $\Sigma_{02}$ blow up to $\sigma={\cal F}^1_{e_2}$.  In this case, we have that $\sigma\cdot S_{01} $ and $\sigma\cdot S_{03}$, are $\pm1$, with opposite signs, so $\sigma\cdot\alpha^+_Z = \pm c_{01}  \mp c_{03} = 0$ as was seen in the proof of Theorem 4.9.

The other possibility is $\ell\subset{\cal I}(f_Z)$, for one of the indeterminate lines in $\Sigma_\gamma$.  This blows up to one of the first fibers $\sigma={\cal F}^1_\zeta$.  In this case, $\sigma$ crosses $\Gamma_1$ transversally, so $\sigma\cdot\Gamma_1=1$.   On the other hand, $\sigma\subset\Gamma_j$ for some $j>1$, so we have $\sigma\cdot\Gamma_j=-1$.   Thus $\sigma\cdot\alpha^+_Z=c'_j-c'_1>0$.
\qed

\noindent {\bf Remark. }  Considering the symmetry between $f$ and $f^{-1}$, we find that ${\cal I}(f^{-1}_Z) = \Sigma_{02}\cup{\cal F}^1_{0\beta\gamma}\cup\bigcup_\zeta{\cal F}^1_\zeta$, where the $\zeta$'s are the intersection points of $\gamma_1$ with the curves $\gamma_2$, $\gamma_3$, $\gamma_7$, and $\gamma_9$.   

If we instead blow up the $\gamma_j$'s in the order $\gamma_1$, $\gamma_2$, \dots, and call the resulting space $\hat Z$.   Then we have ${\cal I}(f_{\hat Z})=\Sigma_{02}\cup{\cal F}^1_{0\beta\gamma}\cup \bigcup_\zeta {\cal F}^1_\zeta$, where the $\zeta\in\gamma_{11}$ are the points of intersection with $\gamma_3$, $\gamma_5$, $\gamma_9$, and $\gamma_{10}$.  Each of these points $\zeta$ is blown up by $f_{\hat Z}$  to a line of the pencil in $\Sigma_C$ passing through $e_3$.

Thus we can apply a similar argument to $\alpha^-_Z$ to obtain the Green current for $f^{-1}_Z$.
\proclaim Corollary 7.8.   There is a positive closed current $T^-_Z$ in the class of $\alpha^-_Z$ with the property: if \ $\Xi^-$\/ is a smooth form which represents $\alpha^-_Z$, then $\lim_{n\to\infty}\delta_1(f)^{-n} f_{Z}^{ -n*}\Xi_Z^- =T^-_Z$  in the weak sense of currents on $Z$.

Next we show what happens to  the invariant fibration when we lift it to $Z$.
Let us set $P_0 = x_0x_1x_2x_3$, and let $P_1$ be a homogeneous quartic polynomial defined in Appendix B.   For $c\in {\bf C}$, let us set  $S_c =\{c P_0 + P_1=0\}$, so the rotor ${\cal R}$ corresponds to $c=\infty$.  Since we have $f(S_c) = S_{\omega c}$, the surface $S_0$ is invariant.

\proclaim Proposition 7.9.  The variety $S_0:=\{P_1=0\}\subset{\bf P}^3$ has singular points at $e_1$, $e_3$ and the fixed points $p_\pm$.  If $p_\pm$ are blown up (in additional to the $e_1$ and $e_3$ which were blown up to construct $Y$), then the strict transform of $S_0$ is a nonsingular $K3$ surface.

\noindent{\it Proof.}   Using the computer, we find that the critical points of $P_1$ occur exactly at $e_1$, $e_3$ and $p_\pm = ( x_\pm, x_\pm, x_\pm) \in {\bf C}^3$ where $x_\pm$ are the roots of $x^2=a + (1+\omega)x$.  ({\sl Mathematica}, for instance, can do this.)  Further, $p_\pm$ are singular points of type $A_1$.  The singular points $e_1$ and $e_3$ are type $A_1$ unless $a=(1+2\omega)/(1-\omega)$, in which case they are type $A_2$.   In either case, it follows (see, for instance,  [EJ, Lemma 3.1 and Remark 3.2]) that $S_0$ is $K3$.  
\qed

%

\proclaim Corollary 7.10. For all but finitely many values of $c\in{\bf C}$, the strict transform of $S_c$ in $Z$ is a nonsingular $K3$ surface. 




Let ${\cal P}\subset Y$ denote the (finite) set of all intersection points of distinct curves $\gamma_j\cap\gamma_k$.   Since the $\gamma_j$ lie in the rotor, we have ${\cal P}\subset{\cal R}$.  The rotor is the union of 8 smooth hypersurfaces which intersect transversally, so the singular locus of ${\cal R}$ is the set where two (or more) of these surfaces intersect.  We will write  ${\cal P}_s$ (resp.\ ${\cal P}_r$) for the points of ${\cal P}$ which are contained in the singular (resp.\ regular) locus of ${\cal R}$.  

While $Z$ itself depends on the order in which the curves $\gamma_j$ are blown up, the following Propositions are valid for any ordering of the blowups.

\proclaim Proposition 7.11.  For  $p\in{\cal P}_r$, there is a unique $c_p\in{\bf C}$ such that $S_{c_p}\subset{  Y}$ is singular at $p$.  This is a conical singularity, and the strict transform $S_{c_p}\subset{  Z}$ contains the first fiber ${\cal F}^1_{c_p}$.

\noindent {\it Proof. }  Without loss of generality, we may choose coordinates $(x,y,z)$ so that $p=0$, $L=z$ near $p$, and ${\cal R}=\{z=0\}$.  Let us suppose that $p\in f_{  Y}^j\Sigma_{BC}\cap f_{  Y}^k\Sigma_{BC}$.   Since the curves $f_{  Y}^j\Sigma_{BC}$ are contained in ${\cal R}$ and intersect transversally, we may suppose that  near $p$ the curves $f_{  Y}^j\Sigma_{BC}$ and $f_{  Y}^k\Sigma_{BC}$ coincide with the $x$- and $y$-axes.  Thus the tangent to $\{M=0\}$ at $p$ is given by $z=0$, so we may suppose that $M=\lambda z+xy+\cdots$.    The surfaces are then $S_c=\{M+cL=0\}=\{ \lambda z+xy + c z+\cdots=0\}$.  The surface $S_c$ is singular if $c=-\lambda$.  We blow up the $x$-axis by the coordinate change  $(x,s,\eta)\mapsto (x,s,s\eta)$.   The first fiber is ${\cal F}^1_p=\{x=s=0\}$.  The strict transforms of the surfaces are $S_c=\{(\lambda+c)\eta +x=0\}$.   The strict transform of the $y$-axis is now the $s$-axis, which is contained in each $S_c$.  Otherwise, the $S_c$'s are disjoint.  The strict transform of $S_{-\lambda}$ contains ${\cal F}^1_p$.  After we blow up the $s$-axis, the surfaces are all disjoint and smooth.  \qed

\proclaim Proposition 7.12.  For $p\in{\cal P}_s$, $S_c$ is smooth at $p$ for all $c\in {\bf C}$.  The first fiber is contained in the rotor:  ${\cal F}^1_p\subset{\cal R}\subset{  Z}$.

\noindent{\it Proof.  }  We may assume that $p$ is a normal crossing of two of the hypersurfaces of ${\cal R}$.  Thus we may choose coordinates $(x,y,z)$ such that $p=0$, and $L=xy$ near $p$.  We may assume that $f_{  Y}^j\Sigma_{BC}$ is the $x$-axis, and $f_{  Y}^k\Sigma_{BC}$ is the $y$-axis.  Since $M$ contains both axes,  we may assume that $M=z+\varphi$, where $\varphi$ is divisible by $xy$.  Thus $S_c=\{M+cL=z+\varphi+cxy=0\}$ is smooth for all $c\in{\bf C}$.    When we blow up the $x$-axis, we use coordinates $(x,s,\eta)\mapsto(x,s,s\eta)$.  The strict transforms are then $S_c=\{ s\eta+ \tilde\varphi+csx=0\}$, where $\tilde\varphi$ is divisible by $xs$.  Dividing this equation by  $s$, we have $S_c=\{\eta+\psi(x,s,\eta)+cx=0\}$, where $\psi(0,s,0)=0$, since $S_c$ contains the $s$-axis (the strict transform of the $y$-axis).   We have ${\cal F}^1_p=\{x=s=0\}$.  Now we blow up the $s$-axis via the coordinates  $(\xi,s,t)\mapsto(\xi  t, s, t)=(x,s,\eta)$.  This gives the new strict transforms $S_c=\{1+\hat\psi(\xi,s,t) + c\xi=0\}$, where $\hat\psi(\xi,s,t)=t^{-1}\psi(\xi t,s,t)$ is regular.  The strict transform of ${\cal F}^1_p$ is now $\{\xi = s= 0\}$, which is disjoint from the $S_c$s.   \qed

If $p'\in\gamma_1\cap\gamma_9$, then there is a unique $c'\in{\bf C}$ be such that $S_{c'}$ is singular at $p'$.  Let $\ell_9'$ denote the line for which $f(\ell_9')=p'$.  By Theorem 7.6 and Proposition 7.11, it follows that $f_{Z}$ maps $\ell_9'$ to the strict transform of $S_{c'}$ inside $Z$.  Thus the total transform of $\ell_9'$ under $f^n_{Z}$ is contained in $S_{\omega^n c'}$.  Let $p''$ denote the other point of $\gamma_1\cap\gamma_9$, and let $c''\in{\bf C}$ denote the corresponding parameter.  Let $\hat S = S_{c'}\cup S_{\omega c'}\cup S_{\omega^2 c'} \cup S_{c''}\cup S_{\omega c''}\cup S_{\omega^2 c''} $ we see that $\hat S$ is a $f_Z$-invariant set which contains $\ell_9'\cup\ell_9''$.  Let ${\cal R}$ denote the strict transform of the rotor in $Z$.  The sets ${\cal R}$, $\hat S$, and $\Sigma_{02}\cup{\cal F}^1_{e_2}$ are totally invariant, and  we break the indeterminacy locus into three sets:
$${\cal I}(f_Z) = (\Sigma_{02}\cup{\cal F}^1_{e_2})\cup ({\cal I}(f_Z)\cap{\cal R})\cup( {\cal I}(f_Z)\cap \hat S)$$
with ${\cal I}(f_Z)\cap {\cal R}=\ell_2\cup\ell_3\cup\ell_7\cup{\cal F}^1_{0\beta\gamma}$, and $ {\cal I}(f_Z)\cap \hat S = \ell_9'\cup\ell_9''$.   We set  
$$\Omega=Z-(\hat S\cup{\cal R}\cup\Sigma_{02}\cup {\cal F}^1_{e_2})$$  
By Propositions 7.11 and 7.12, ${\cal R}$ is disjoint from the strict transform of each $S_c$.  Thus $f_Z$ is regular on $\Omega$, and $\Omega$  is invariant under $f_Z$.

\proclaim Proposition 7.13.  For every $S_c$ in $\Omega$ the  dynamical degree of the restriction is $\delta(f^3_c)=\delta_1(f)^3$.  

\noindent{\it Proof.}  Let us denote $\Gamma$ a hypersurface in ${  Z}$ whose cohomology class in $H^{1,1}({  Z}) $ is  $H_{  Z}$. It follows that the degree of $f_Z^{-3n} \Gamma$ grows like $\delta_1(f)^{3n}$. On the other hand $S_c\subset \Omega $ does not contain an irreducible component of the indeterminacy locus for $f_{  Z}$. It follows that we have $S_c \cap (f_Z^3)^{-n} \Gamma = (f_{  Z}^3)^{-n} ( S_c \cap \Gamma)$. Because $S_c$ is non-singular and $f_Z^3$ is pseudo-automorphism, the degree of $(f_{  Z}^3)^{-n} ( S_c \cap \Gamma) = (f^3_c)^{-n}( S_c \cap \Gamma)$ is $4\, \delta_1(f)^{3n}$. Thus the dynamical degree of $f^3_c$ is $\delta_1(f)^3$.\qed

Using the fact that $f_Z$ is regular on the large invariant set $\Omega$, we avoid the difficulties that can occur in defining the entropy of a map (see [G1]).
 \proclaim Theorem 7.14.  The entropy of $f$ is $\log\delta_1(f)$.

\noindent{\it Proof.}   Since $f$ is equivalent to a pseudo automorphism,  and $f^*_Z$ is conjugate to $(f^*_Z)^{-1}$, both the first and the second dynamical degrees are equal. Combining the result in [DS] and the fact that $h_{top}(f) \ge h_{top}(f_0)$, we have the inequality $$ \log \delta_1(f) \ge h_{top}(f) \ge h_{top}(f_0) = \log \delta_1(f)$$ 
which gives the result.  \qed

 Since $\Xi^\pm$ and $f_Z$ are regular on $\Omega$,  the potential of $g^\pm$ is continuous on $\Omega$.  Thus may define the wedge product $T_2:=T^+\wedge T^-$ as a positive, closed (2,2)-current on $\Omega$, and we have:    
\proclaim Proposition 7.15.   $\lim_{n\to\infty}\delta_1(f)^{-2n} {f_{Z}^*}^{n}\Xi^+ \wedge {f_{Z}^*}^{-n}\Xi^-= T_2$ exists as a (2,2)-current on $\Omega$.
%

We have seen that the restrictions $f^3|_{S_c}$ are automorphisms, and there are invariant currents $\mu^\pm_c$ on $S_c$, as well as invariant measures $\mu_c:=\mu^+_c\wedge\mu^-_c$ (see [C]).  The following property leads us to consider $T^\pm$ and $T_2$ as  the ``bifurcation currents'' for the family $\{f^3|_{S_c}\}$ (see [DuF]).

\proclaim Theorem 7.16.  For $S_c\subset\Omega$ the slices by $S_c$ are well-defined and give the corresponding dynamical objects:  $T^\pm|_{S_c}=\mu^\pm_c$, and $T_2|_{S_c}=\mu_c$.

\noindent{\it Proof. }  If we set $h=f^3$, then the class $[S_c]$ is invariant under $h^*$.   Thus $\alpha^+\cdot[S_c]\in H^{1,1}(S_c)$ is a class that is expanded by a factor of $\delta_1(f)$.  It follows that the restriction $\Xi^+|_{S_c}$ gives the expanded class, and this converges to $\mu^+_c$.  Similarly, the normalized pullbacks/push-forwards of $\ \Xi^+\wedge\Xi^-$ on $S_c$ will converge to $\mu_c$.  \qed

\proclaim Theorem 7.17.  For generic $c',c''$, the maps $f^3|_{S_{c'}}$ and $f^3|_{S_{c''}}$ are not smoothly conjugate, and the surfaces $S_{c'}$ and $S_{c''}$  are not isomorphic.

\noindent{\it Proof. }  There is an invariant 6-cycle of curves, $\Gamma_j$, $j=0,\dots,5$ for $f$.  For generic $c$,  $\Gamma_j\cap S_c$ is a saddle 2-cycle for $f^3|_{S_c}$.  The multipliers of this saddle cycle are not constant in $c$, so the maps $f^3|_{S_c}$ are not smoothly conjugate.  Since the automorphism group of $S_c$ is disconnected, we see that the family $\{S_c\}$ cannot consist of surfaces which are all isomorphic to each other.  \qed

\noindent{\bf Remark.}  If $a_2\ne 1$ is a primitive 5th root of unity and $a_0= b_0=0$, then we may repeat most of the arguments in this section for this map.  In particular, we have:

\proclaim Theorem 7.18.  If $a_2$ is a primitive 5th root of unity, and $a_0=b_0=0$, then $f$ is equivalent to a pseudo-automorphism, and the dynamical degrees $\delta_1(f)=\delta_2(f)\approx 1.3211018 > 1$ are the largest root of $t^{19}(t^3-t-1)+t^3+t^2-1$.  The entropy of $f$ is $\log \delta_1(f) >0$. Furthermore there are two quartic polynomials which are invariant in the sense of (B.1).  This gives a family of $K3$ surfaces which are invariant under $f^5$.   

\bigskip\noindent{\bf\S 8. Pseudo-automorphisms which are completely integrable. }  Let us consider two cases for maps of the form (2.2):
$$ \alpha=(a,0,1,1),\  a\ne1,\  {\rm\ and \ } \beta=(0,1,0,0)\eqno(8.1a)$$
$$\alpha=(0,0,\omega,1), \ \omega^3=1,\  \omega\ne1,\ {\rm\ and\ } \beta=(0,1,0,0)\eqno(8.1b)$$
The map (8.1a) has been extensively studied under the name Lyness process.    The maps (8.1a) and (8.1b) exhibit similarities to the maps in the previous section:  they are critical maps, and the iterates of the  critical image $\Sigma_{BC}$ go ``once around'' the rotor and land on $\Sigma_{\beta\gamma}$.  The difference with \S7 is that $f_{  Y}^4\Sigma_{BC}={\cal F}_{0\beta\gamma}$ is an indeterminate curve, and by Lemma 4.4 this fiber is mapped to ${\cal F}_{0BC}$, that is, $f^4_{  Y}\Sigma_{BC} = {\cal F}_{0\beta\gamma} \subset S_{01}$ and $f^5_{  Y}\Sigma_{BC} = {\cal F}_{0BC} \subset S_{03}$. Thus $\Sigma_{BC}$ arrives at $\Sigma_{\beta\gamma}$ one step faster than was the case in \S7.

Let $\pi:{  Z}\to {  Y}$ denote the space obtained by blowing up the orbit $f^j\Sigma_{BC}$, $f^{-j}\Sigma_{\beta\gamma}$, $0\le j\le 4$ (one curve less than the construction in \S7).

\proclaim Theorem 8.1.  The induced map $f_{  Z}$ is a pseudo-automorphism, and the iterates of $f$ have quadratic degree growth.

\noindent{\it Proof. }  Since $f^9_{  Y} \Sigma_{BC} = \Sigma_{\beta\gamma}$ and $f_{  Y}^4 \Sigma_{BC} = {\cal F}_{0\beta\gamma},  f_{  Y}^5 \Sigma_{BC} = {\cal F}_{0BC}$, we see that $f_{  Y}$ satisfies the condition in Theorem 5.1. This Theorem then follows from Theorems 5.1 and  5.3. \qed

\proclaim Proposition 8.2.  In cases (8.1a) and (8.1b), the induced rotor map $f^8|_{\Sigma_3}$ has linear degree growth.  This map is not birationally conjugate to a surface automorphism. 

\noindent{\it Proof.} 
In the case  (8.1b), the restriction of $f_{  Y}^8$ to $\Sigma_3$ is given by setting $a_0=0$ in (7.3), so we find the degree $2$ birational map : 
$$f_{  Y}^8|_{ \Sigma_3}: [x_0:x_1:x_2:0] \mapsto [ x_0 \omega^2 x_2 : x_1 ( x_1+\omega^2 x_0 + \omega^2 x_2)
:\omega^2 x_2  (x_1+ \omega^2 x_2) : 0] $$
This map has three distinct exceptional lines. Two of exceptional lines are mapped to fixed points $[1:0:-1:0]$ and $[0:1:-a:0]$. The remaining exceptional line is mapped to a point of indeterminacy $e=[1:-1:0:0]$. We let ${  W}$ be the blowup space obtained by blowing up $\Sigma_3$ at $e$. The induced map has only two exceptional lines which are mapped to fixed points and therefore the induced map is algebraically stable. The action on $Pic$ is given by the matrix $\pmatrix{2&1\cr -1& 0} $ which has an eigenvalue $1$ with $2\times 2$ jordan block. It follows that the degree of restriction map grows linearly.   

The analysis in the case (8.1a) is essentially the same.  The induced rotor map is now:
$$f_{  Y}^8|_{\Sigma_3}:  [x_0:x_1:x_2:0]\mapsto [x_0(ax_0 + ax_1 + x_2): x_1(x_0+x_1+x_2): x_2(ax_0 + x_1 + x_2): 0]$$
This map has three exceptional lines.  Two of them are mapped to fixed points $[1:0:-1]$ and $[0:1:-a]$.  The third exceptional line is mapped to  $[1:-1:0]$, which is indeterminate. After we blow up the point $[1:-1:0]$, the induced map is algebraically stable and the action on $Pic$ has an eigenvalue $1$ with $2\times 2$ jordan block. 

Finally, since the restriction of $f_{  W}$ to the rotor has linear degree growth.  It follows from [DiF] that this restriction is not an automorphism.   \qed

We consider first the Lyness map, i.e.,  case (8.1a).  This is known to be integrable, and the invariant polynomials are given in [CGM1] and [KoL].   These invariant polynomials, which satisfy (B.1) with $t=1$, are:
$$\eqalign{&Q_0 = x_0x_1x_2x_3\cr &Q_1= (a x_0 + x_1+x_2+x_3) (x_0+x_1) (x_0+x_2) (x_0+x_3)\cr &Q_2 =(x_0 ( a x_0+x_1+x_2+x_3)+x_1 x_3) (x_0+x_1+x_2) (x_0+x_2+x_3).\cr}  \eqno(8.2a)$$
The set $\{Q_0=0\}$ gives an invariant 8-cycle of rational surfaces, which is the rotor  ${\cal R}\subset {  Y}$.  (Although $Q_0=0$ consists of 4 irreducible components in ${\bf P}^3$, it yields an 8-cycle inside $Y$ because these components map through the indeterminacy locus, which is blown up to yield an additional 4 divisors.)  The set $\{Q_1=0\}$ gives an invariant 4-cycle, and $\{Q_2=0\}$ gives an invariant 3-cycle; the components of the 8-, 4-, and 3-cycles are rational surfaces.   As we observed in \S4,  $f_{  Y}$ induces dominant maps on each of these cycles.   And as in Proposition 8.2, we may show that the restriction of $f^4$ to the 4-cycle, and the restriction of $f^3$ to the 3-cycle both have linear degree growth.


Let us define the surfaces $S_c=\{Q_c=0\}$ with $Q_c:=c_0Q_0+c_1Q_1+c_2Q_2$.  If we also write $S_c$ for its strict transform inside ${  Z}$,  we have $fS_c=S_c$ 

\proclaim Theorem 8.3.  For generic $c$, the surface $S_c$ is an irreducible $K3$ surface.

\noindent{\it Proof. }  For generic $c$, we find that $S_c$ has $16$ singular points:  two of them are $e_1,e_3$, which are type $A_2$,  and there are $14$ more which are of type $A_1$.  In the construction of $Z$, we blew up $e_1$ and $e_3$.  Then we blew up $f^j\Sigma_{BC}$, $0\le j\le 10$, and the other 14 singular points are contained in these curves.  It follows that the strict transform of $S_c$ inside $Z$ is smooth and thus $K3$.   \qed

\proclaim Theorem 8.4.  For generic $c$ and $c'$, the intersection $S_c\cap S_{c'}$ is an elliptic curve.  The restriction of $f^3$ to $S_c$ has quadratic degree growth.

\noindent{\it Proof.}  Since $S_c$ is a $K3$ surface, it has trivial canonical bundle.  Thus the birational map $f^3$ of $S_c$ must be an automorphism.  For generic $c$ and $c'\ne c$, the intersections $S_c\cap S_{c'}$ give an invariant fibration of $S_c$. Since $f^3|S_c$ is an automorphism, then by [DiF] the intersection  $S_c\cap S_{c'}$  is an elliptic curve and the restriction of $f$ to the family of $K3$ surfaces has quadratic degree growth.  \qed

The map (8.1b) is similar.  In this case the solutions to (B.1) take the form:
$$\eqalign{ & R_0 = x_0x_1x_2x_3 \cr
& R_1=(x_0+\omega x_1)(x_0+\omega x_2)(x_0 + \omega x_3)(x_1+\omega^2 x_2+\omega x_3)\cr
& R_2 =\omega x_1x_3 (x_0+\omega x_1)(x_0+\omega x_3)\cr &\phantom{R_2=AA} + \omega^2 x_0x_2  \left( x_0(x_1+\omega x_3) + x_2 (\omega x_1 + x_3) + \omega^2 x_0 x_2\right)\cr} \eqno(8.3)$$
where $t_{R_0}=1$, $t_{R_1} =\omega^2$, and $t_{R_2}=\omega^2$.  As before, we see that  $f_{  Z}$ will have an invariant 8-cycle given by the rotor ${\cal R}\subset{  Z}$.  And $\{R_1=0\}$ will give a 4-cycle of rational surfaces.  For generic $c$, the singularities of the surface $S_c=\{\sum c_j R_j=0\}$ are $e_1,e_3$ (type $A_2$) and $e_2$ (type $A_1$).  As in Theorems 8.3 and 8.4, we have:

\proclaim Theorem 8.5.  In case (8.1b): for generic $c$, $S_c$ is a $K3$ surface,  $f^3$ is an automorphism of $S_c$ with quadratic growth, and the intersections $S_c\cap S_{c'}$ are elliptic curves.

%

\bigskip\noindent{\bf\S A.  Appendix: Computing the Characteristic Polynomial for $f_{  Z}^*$. }  Let us consider a critical map $f$ satisfying condition (5.1), and let $m_u$, $m_d$, $m_s$, $d_j$, $u_j$ and $N$ be the numbers defined in \S5.   We define the $(N+5)\times(N+5)$ matrix
$$\pmatrix { 2 & 0&1&0&1&0 &  \cdots &  0&1 \cr 
                                   -1 & 0&-1&0&0&0& \cdots & 0&-1 \cr 
                                   0&1&-1&0&0&0& \cdots& 0&0\cr
                                   -1&0&-1&0&-1&0& \cdots & 0&0\cr
                                   -1&0&-1&1&-1&0&\cdots &  0&0\cr
                                   -1 & 0 &0&0&-1&0&\cdots & 0&-1\cr
                                  * &  & & & & & & &  *\cr
                                    \vdots &  & & & &   & & & \vdots\cr
                                    *&  & & & & & &   &*\cr} \eqno(A.1)$$
 where the $*$'s indicate that the $7^{\rm th}$ through the $N+5^{\rm th}$ rows remain to be specified.   We will define the $j^{\rm th}$ row $r_j$ in terms of the elements $e_k$, which are vectors of length $N+5$ in which the $k^{\rm th}$ entry is 1, and all other entries are 0:
 \item{(a)} if $j=N- d_i$ for some $i = 1, \dots, m_d$, then  $r_{j+6} = e_{j+5} - e_{N+5}$
 \item{(b)} if $j=N- u_i$ for some $i = 1, \dots, m_u$, then $r_{j+6} = -e_1 - e_5 + e_{j+5} - e_{N+5}$                         
\item{(c)} if $j= N-m_s$, 
$$\eqalign{ r_{j+2} = -e_1 -e_3 + e_{j+1} -e_{N+5},\ \ \ &  r_{j+3} = -e_3 +e_{j+2} \cr
r_{j+4} = -e_1-e_3 - e_5 + e_{j+3} - e_{N+5}, \ \ \ &r_{j+5} = -e_1 - e_3 - e_5 + e_{j+4} \cr
r_{j+6} = -e_5 + e_{j+5} \ \ \ & \cr }  $$ 
\item{(d)} Otherwise, $r_{j+6} = e_{j+5}$.

Let  $\pi:{  Z}\to{\bf P}^3$ be the space constructed in Theorems 5.1 and 5.3.

\proclaim Proposition A.1.  The matrix $(A.1)$  represents $f_{  Z}^*$.

\noindent{\it Proof. } There are three cases to consider.  Although we define a different basis in each case, the matrix $(A.1)$ representing $f^*_{  Z}$ is the same.

\noindent Case (I) :  There is no $1 \le j\le N$ such that $f^j_{  Y} \Sigma_\gamma \subset {\cal F}_{0\beta\gamma} \cup \Sigma_{\beta\gamma}$. \hfill \break In this case for all $j \ne d_i, u_k, m_s, m_{s+1}, \dots, m_{s+4}$ , $ 1 \le i \le m_d, 1 \le k \le m_u$ we have $ \pi {\cal F}_j \not\subset \Sigma_0 \cup \Sigma_\beta\cup \Sigma_\gamma$ and therefore
$$\eqalign{ f^*_{  Z} H_{  Z}  & = 2 {\cal H}_{  Z} - E_1 - S_{01}- E_3  - \sum_{i=1}^{m_u} {\cal F}_{u_i} - {\cal F}_{m_s+1}  - {\cal F}_{m_s+2}  - {\cal F}_{m_s+4}\cr
\{\Sigma_0\} &= {\cal H}_{  Z} - E_1 -S_{03} - S_{01}- E_3  - {\cal F}_{m_s+1}  - {\cal F}_{m_s+2} - {\cal F}_{m_s+3} - {\cal F}_{m_s+4}\cr
 \{\Sigma_\beta\} &= {\cal H}_{  Z}  - S_{01}- E_3 - {\cal F}_N  - \sum_{i=1}^{m_u} {\cal F}_{u_i} - {\cal F}_{m_s}  - {\cal F}_{m_s+1} - {\cal F}_{m_s+2} \cr 
 \{\Sigma_\gamma\} &= {\cal H}_{  Z} - E_1 -{\cal F}_N  - {\cal F}_{m_s+1}  - \sum_{i=1}^{m_u} {\cal F}_{u_i}  - \sum_{i=1}^{m_d} {\cal F}_{d_i}. \cr}$$
Since we have $f^*_{  Z} : E_1 \mapsto S_{03} \mapsto \{ \Sigma_0\}, \ \ S_{01} \mapsto E_3 \mapsto \{ \Sigma_\beta\}, $ $ {\cal F}_j \mapsto {\cal F}_{j-1}$ for all $j = 2, \dots, N$, and ${\cal F}_1 \mapsto \{\Sigma_\gamma\}$ using the ordered basis  $ \{ H_{  Z}, E_1, S_{03}, S_{01},E_3, {\cal F}_N, {\cal F}_{N-1}, \dots, {\cal F}_2, {\cal F}_1\}$  for ${\rm Pic}({  Z})$ we see that $(A.1)$ is the matrix representation for $f_{  Z}^*$.

\smallskip
\noindent Case (II) : There are $\kappa$ positive integers $ 1 < s_1 < \cdots < s_\kappa < N$ such that  for $j = 1, \dots, \kappa$ $f^{s_j} _{  Y} \Sigma_\gamma \subset \Sigma_\beta \setminus \ell_\beta \cup \Sigma_{\beta\gamma} \cup \ell_\beta' $ where $f^3_{  Y} \ell_\beta' = E_3 \cap \Sigma_0$. \hfill\break
For this case let us use the ordered basis $ \tilde{\cal B} = \{ H_{  Z}, \tilde E_1, \tilde S_{03}, \tilde S_{01}, \tilde E_3, {\cal F}_N, {\cal F}_{N-1}, \dots, {\cal F}_2, {\cal F}_1\}$  for ${\rm Pic}({  Z})$ where $\tilde E_3=  E_3 + \sum_{i=1}^{\kappa} {\cal F}_{s_i + 1}$, $\tilde S_{01}=  S_{01} + \sum_{i=1}^{\kappa} {\cal F}_{s_i + 2}$,$\tilde S_{03}=  S_{03} + \sum_{i=1}^{\kappa} {\cal F}_{s_i + 4}$, and $\tilde E_1=  E_1 + \sum_{i=1}^{\kappa} {\cal F}_{s_i + 5}$. Using this new ordered basis we can see that 
$$\eqalign{f^*_{  Z} : \tilde E_1 \mapsto \tilde S_{03}& \mapsto  \{ \Sigma_0\} + \sum_{i=1}^\kappa {\cal F}_{s_i+3}\cr & = {\cal H}_{  Z} - \tilde E_1 -\tilde S_{03} - \tilde S_{01}- \tilde E_3  - {\cal F}_{m_s+1}  - {\cal F}_{m_s+2} - {\cal F}_{m_s+3} - {\cal F}_{m_s+4}\cr}$$
In a similar way we may compute  $ f^*_{  Z}$ of $H_{  Z}, \tilde S_{01}$, $\tilde E_3$ and ${\cal F}_N$ and see that the matrix representation with  $ \tilde{\cal B} $ is given by $(A.1)$.

\smallskip
\noindent Case (III) : There are $\tau$ positive integers $ 1 < q_1 < \cdots < q_\tau< N$ such that $f^{q_j} _{  Y} \Sigma_\gamma \subset  \ell_\beta' $ for $j = 1, \dots, \tau$. \hfill\break
Let us consider the ordered basis $ \hat{\cal B} = \{ H_{  Z}, \hat E_1, \hat S_{03}, \hat S_{01}, \hat E_3, {\cal F}_N, {\cal F}_{N-1}, \dots, {\cal F}_2, {\cal F}_1\}$  for ${\rm Pic}({  Z})$ where $\hat E_3=  \tilde E_3 + \sum_{i=1}^{\tau}( {\cal F}_{q_i + 1} + {\cal F}_{q_i + 3})$, $\hat S_{01}=  \tilde S_{01} + \sum_{i=1}^{\tau}( {\cal F}_{q_i + 2} + {\cal F}_{q_i + 4})$, $\hat S_{03}=  \tilde S_{03} + \sum_{i=1}^{\tau} ({\cal F}_{q_i + 4} +{\cal F}_{q_i + 6})$, and $\hat E_1=  \tilde E_1 + \sum_{i=1}^{\tau} ({\cal F}_{q_i + 5} +{\cal F}_{q_i + 7})$. Since $f_Y^2 \ell'_\beta = \Sigma_\beta \cap S_{01}$, we have 
$$\eqalign{ \{\Sigma_\beta\}  &=  {\cal H}_{  Z} -  \tilde S_{01}- \tilde E_3  - \sum_{i=1}^\tau ( {\cal F}_{q_i} + {\cal F}_{q_i+1}+ 2 {\cal F}_{q_i+2}+ {\cal F}_{q_i+3}+ {\cal F}_{q_i+4} ) \cr
& \qquad\qquad -{\cal F}_N- \sum_{i=1}^{m_u} {\cal F}_{u_i} - {\cal F}_{m_s}  - {\cal F}_{m_s+1} - {\cal F}_{m_s+2}\cr
& = {\cal H}_{  Z} -  \hat S_{01}- \hat E_3  - \sum_{i=1}^\tau ( {\cal F}_{q_i} +  {\cal F}_{q_i+2} )-{\cal F}_N- \sum_{i=1}^{m_u} {\cal F}_{u_i} - {\cal F}_{m_s}  - {\cal F}_{m_s+1} - {\cal F}_{m_s+2}.}$$
It follows that we have 
$$\eqalign{ f^*_{  Z} : \hat S_{01} \mapsto \hat E_3 \mapsto &   \{ \Sigma_\beta\} +  \sum_{i=1}^\tau ( {\cal F}_{q_i} +  {\cal F}_{q_i+2} )\cr
& =  {\cal H}_{  Z}  - \hat S_{01}- \hat E_3 - {\cal F}_N  - \sum_{i=1}^{m_u} {\cal F}_{u_i} - {\cal F}_{m_s}  - {\cal F}_{m_s+1} - {\cal F}_{m_s+2} \cr }$$
For the other basis elements, computations are essentially identical and thus we see that $(A.1)$ represents $f^*_{  Z}$ with respect to the ordered basis $\hat {\cal B}$. 
\qed

According to the previous Proposition, we see that the characteristic polynomial of $f^*_{  Z}$ only depends on  $m_u$, $m_d$, $m_s$, $d_j$, $u_j$ and $N$.
 
\proclaim Lemma A.2. The characteristic polynomial of $f^*_{  Z}$ is given by $$\pm t^{N-1}(t^2+1) \left[ (Q_1 - Q_4) t^3  +( 2 Q_1 - Q_2 - Q_3 -Q_4) t^2 + ( Q_1 - Q_3) t + Q_4 \right].$$

\noindent{\it Proof.}
We subtract $tI$ from the matrix $(A.1)$ and perform a sequence of row operations on it.  Step  (i): we  add or subtract the $6^{\rm th}$ row to the rows whose last entry is $1$ or $-1$ and then (ii)  for $j=1,\dots,N-1$, we subtract $1/t^j$ times the  $N+4-j^{\rm th}$ row  from $6^{\rm th}$ row. This gives
 $$ {\rm det} (f_{  Z}^* - t \,I) = {\rm det} \pmatrix {A&0\cr *&B}$$
 where  
 $$A= \pmatrix {1-t&0&1&0&0&-t\cr 1-t&-t&0&0&1&0\cr 0&1&-1-t&0&0&0\cr -1&0&-1&-t&-1&0\cr -1&0&-1&1&-1-t&0\cr Q_1 & 0 &Q_2 &0&Q_3 & Q_4\cr},\ \ B= \pmatrix{ -t& 0 & 0&\cdots &0& 0\cr 1&-t &0& \cdots &0&0 \cr 0 &1& -t & \cdots &0&0\cr \vdots & & \ddots & \ddots & & 0\cr 0 & & & \ddots& -t&0\cr 0 & & & & 1& -t}$$
with $Q_1, Q_2, Q_3, Q_4$ as in \S 5.  We have $ {\rm det} (f_{  Z}^* - t \,{\rm Id}) =(-1)^{N-1} t^{N-1} {\rm det}(A)$, and we evaluate ${\rm det}(A)$ to obtain the polynomial given above.   \qed
\bigskip

\noindent{\bf \S B.  Appendix:  Invariant Polynomials.  }   We will look for polynomials $P(x) = \sum a_I x^I$ which are invariant in the sense that 
$$P\circ f= t\cdot j_f\cdot P \eqno(B.1)$$
where $t\ne0$ is constant, and $j_f= 2x_0 (\gamma\cdot x)(\beta\cdot x)^2$ is the Jacobian determinant.   If $P$ and $Q$ are solutions to (B.1) with multipliers $t_P$ and $t_Q$, then $\varphi=P/Q$ is a rational function with the invariance property: $\varphi\circ f= (t_P \,t_Q^{-1})\, \varphi$.  If  $P$ is a solution to (B.1), then $P$ defines a meromorphic 3-form $\Omega_P$: on the set $x_0\ne0$, it is given by $P(1,x_1,x_2,x_3)^{-1} dx_1\wedge dx_2\wedge dx_3$.  This is invariant in the sense that $f^*\Omega_P=t_P^{-1}\Omega_P$.  It follows that $\{P=0\}$ is an $f$-invariant surface which represents the canonical class in ${\bf P}^3$ and its strict transforms are invariant surfaces which represent the canonical classes in ${  Y}$ and ${  Z}$.


The equation (B.1) can be rewritten as a system of linear equations for the coefficients of the monomials in $P$.  This system can be solved directly for all the maps in \S7 and \S8.  For instance, in \S7  $\omega$ is a non-real root of unity and $a_0=a\ne0$, and  we find a solution for $t=\omega^2$:
%


$$\eqalign{
P_1& =(1-\omega)\left(a^2 x_0^4  +(1 + a) x_0 x_1 x_2^2 
+ x_1^2 x_3^2  +a x_1 x_2 x_3^2 \right ) 
  \cr 
 &-(2 + \omega )\left( x_0 x_2^3  + (1 + a) x_0 x_1^2 x_3  +
 a x_1 x_2^2 x_3  +a x_0^2 x_3^2 \right) 
  \cr
 &
 +(1 + 2 \omega ) \left( a x_0^2 x_1^2 
 + a x_0 x_1^2 x_2  + 
 a x_1^2 x_2 x_3 + a x_0 x_2 x_3^2 \right) 
+ a x_0^3 x_1 (1 + a + 2 \omega  - a \omega ) 
 \cr
 &+ (1- 2 a + 2 \omega  - a \omega )\left( (1 +  a) x_0^2 x_1 x_3 +
 x_0 x_2^2 x_3\right ) + 
 x_0^2 x_2^2 (1 - a + 2 \omega  + a \omega )  \cr
 &-  (2 - a + \omega  + a \omega ) \left((1 +    a) x_0^2 x_1 x_2
 + x_0 x_1 x_3^2 \right )  
 +a x_0^3 x_3 (1 - 2 a - \omega  - a \omega ) \cr
 &
 + (1 +     a) x_0^2 x_2 x_3 (1 + a - \omega  + 2 a \omega ) + 
 a x_0^3 x_2 (2 + a + \omega  + 2 a \omega ) \cr
}$$

\medskip

\noindent{\bf \S C. Appendix: The Rotor Map.} Let $g:=f^8_Z|_{\Sigma_3}$  denote the rotor map restricted to $\Sigma_3$, which is written in coordinates in (7.3).  
%
By factoring the jacobian determinant, we see that there are four exceptional curves. 
$$\eqalign{ & {\cal C}_1  = \{ a x_0 + \omega x_2 =0\} \cr
&{\cal C}_2 = \{ a x_0 + a x_1 + \omega x_2 =0\} \cr 
&{\cal C}_3 = \{ a \omega x_0 + x_1 + \omega^2 x_2 =0\}\cr
&{\cal C}_4= \{ a \omega x_0^2 + a \omega x_0 x_1 + a \omega x_0 x_2 + \omega^2 x_0 x_2 + x_1x_2 + \omega^2 x_2^2 =0\}}$$

\proclaim Lemma C.1. If $ a \ne \omega^j$ and $a^j  \ne \omega^{j\pm2}$ for all $j \ge 2$ then $g$ is not birationally conjugate to an automorphism.

\noindent{\it Proof. } The exceptional curves ${\cal C}_2$ and ${ \cal C}_4$ mapped to a three cycle : $g: {\cal C}_2 \mapsto [0:1: -a \omega] \mapsto [0:1: - a] \mapsto [0:1: - a \omega^2] \mapsto [0:1: - a \omega]$ and $g: {\cal C}_4 \mapsto [1:0 - \omega^2] \mapsto [1:0 :-\omega] \mapsto [1:0 : -1] \mapsto [1:0: - \omega^2]$.  For ${\cal C}_3$ we see that $g^j {\cal C}_3 = [1: - \omega^2 (\omega/a)^{j-1}:0]$ for all $j \ge 1$. It follows that these three curves have orbits that do not  encounter the indeterminacy locus of $g$. The remaining exceptional curve ${\cal C}_1$ mapped to $e_1= [0:1:0]$, which is indeterminate. We let $W$ be the space obtained by blowing up $\Sigma_3$ at $e_1$, and we let $E_1$ be the corresponding exceptional divisor. Under the induced map $g_W$ we have $g_W (E_1) = E_1$ and the orbit of the strict transform of ${\cal C}_1$ remains in $E_1$ and does not encounter the indeterminacy locus of $g_W$. 

Now if $H$ denote the class of a generic line in $W$, then $\langle H, E_1\rangle$ is an ordered basis for $Pic(W)$. The action on $Pic$ is given by the matrix $g^*_W = \pmatrix { 3 & 1 \cr -1 & 0}$. The largest eigenvalue is $ \lambda = (3 + \sqrt{5})/2$ and invariant class is given by $ \theta  = \lambda H - E_1$. Since $ \theta^2 = \lambda^2 -1 \ne 0$, it follows from [DF, Theorem 5.4] that $g$ is not birationally conjugate to an automorphism. \qed

\proclaim Lemma C.2. If $ a^j = \omega^{j-2}$ for some $j \ge 2$ then $g$ is not birationally conjugate to an automorphism.

\noindent{\it Proof. } In case $a^j = \omega^{j-2}$ for some $j \ge 2$, the orbits of three exceptional curves ${\cal C}_2, {\cal C}_3,$ and ${\cal C}_4$ are the same as the previous Lemma. After we blowup $e_1$ on $\Sigma_3$, the strict transformation of ${\cal C}_1$ mapped to a point of indeterminacy after $j$-th iteration of $g_W$. We let $W_2$ be the space obtained by blowing up $W$ at $g_W^k {\cal C}_1$ for $k=1, \dots, j$ and we let $F_k$, $1\le k \le j$ be the corresponding exceptional divisors. Under the induced map $g_{W_2}$, the exceptional line ${\cal C}_1$ is removed and the orbits of remaining three exceptional curves do not encounter the indeterminacy locus of $g_{W_2}$.  

Let $\langle H, F_j, F_{j-1}, \dots, F_1, E_1\rangle $ be the ordered basis for $Pic(W_2)$. The characteristic polynomial of the action on $Pic$ is given by $t^{j+2} - 4 t^{j+1} + 3 t^j+t^2 -2 t +1$. It follows that the dynamical degree is not a Salem number. Thus by [DiF], $g$ is not birationally conjugate to an automorphism. \qed

\proclaim Lemma C.3. If $ a^j = \omega^{j+2}$ for some $j \ge 2$ then $g$ is not birationally conjugate to an automorphism.

\noindent{\it Proof. } When $a^j = \omega^{j+2}$, the orbit of ${\cal C}_3$ is different from Lemma C.1, that is $g^{j+1} {\cal C}_3 = [1:-1:0]$, which is indeterminate. We let $W_3$ be the space obtained by blowing up $\Sigma_3$ at $e_1$ and $g^k {\cal C}_3, 1 \le k \le j+1$, and we let $E_1$ and $F_k, 1 \le k \le j+1$ be the corresponding exceptional divisors. Using the ordered basis $\langle H, F_{j+1}, F_{j}, \dots, F_1, E_1\rangle $ for $Pic(W_3)$, we see that the characteristic polynomial of the action on $Pic$ is given by $t^{j+3} - 3 t^{j+2}+ t^{j+1} +t $. Similarly as in Lemma C.2, the dynamical degree is not a Salem number and therefore $g$ is not birationally conjugate to an automorphism.\qed

\proclaim Lemma C.4. If $a = \omega$ then $g$ is not birationally conjugate to an automorphism.

\noindent{\it Proof. } In this case we see that ${\cal C}_2$ is mapped to a point of indeterminacy under $2$ iterations and ${\cal C}_4$ is also mapped to a point of indeterminacy under $3$ iterations. After we blowup $e_1$, we can check that the orbits of other two remaining exceptional lines does not encounter the indeterminacy locus. After we blow up the orbit of ${\cal C}_2$ and the orbit of ${\cal C}_4$, we see that the dynamical degree of $g$ is given by the largest root of the polynomial $ t^3  - t^2 -2t-1$. Again since this number is not a Salem number we have our result. \qed

\proclaim Lemma C.5. If $a = \omega^2$ then $g$ is not birationally conjugate to an automorphism.

\noindent{ \it Proof. } If $a = \omega^2$ the each component of $g$ has the same factor $ x_0 + x_1 + \omega^2 x_2$. It follows that the restriction of $f^8_{  Y}$ to $\Sigma_3$ is a degree $2$ birational map. There are two exceptional lines and both exceptional lines are mapped to points of indeterminacy. After we blowup the points on the orbits of three exceptional lines, we see that the induced map has one exceptional line which is mapped to a point of indeterminacy. Once we blow up this point of indeterminacy, we see that the induced map has no exceptional lines and therefore the induced map is algebraically stable. Furthermore the characteristic polynomial of the action on $Pic$ is $t(1+t)(t-1)^3$ and the action on $Pic$ has $2\times2$ Jordan block. It follows that the degree of $g$ grows linearly. According to [DiF], we have that  $g$ is not birationally conjugate to an automorphism.\qed

\proclaim Lemma C.6. If $a=1$ then the degree $g$ grows quadratically.

\noindent{ \it Proof.} For this case all four exceptional curves are mapped to points of indeterminacy: $g:  {\cal C}_1 \mapsto e_1, \ {\cal C}_2 \mapsto [0:1: -\omega], {\cal C}_3 \mapsto [1:-\omega^2:0] \mapsto [1:-1:0]$ and $g: {\cal C}_4 \mapsto [1:0: -\omega^2]$. We let $Z$ be the space obtained by blowing up $\Sigma_3$ at all five points in the orbit of exceptional curves and we let $E_1$, $Q_2,Q_3,Q_4$, and $Q_5$ be the corresponding exceptional divisors. Under the induced map $g_Z$, there is a unique exceptional line which is the strict transformation of ${\cal C}_1$. We see that $g_Z {\cal C}_1$ is a point of indeterminacy of $g_Z$. By blowing up one more point on $E_1$ we make the induced map an algebraically stable. Let us denote $Q_1$ the exceptional divisors corresponding to the point blow ups on $E_1$.  Let us use $\langle H, F_1,F_2,F_3,F_4,F_5,E_1\rangle$ as the ordered basis of $Pic$. The characteristic polynomial of the action on $Pic$ is given by $(t-1)^4 (t+1) (t^2+t+1)$ and the matrix representation of the action on $Pic$ has $3 \times 3$ Jordan block. It follows that the degree of $g$ grows quadratically. \qed

\bigskip

\centerline{\bf References}
\medskip

\item{[Ba]}  T. Bayraktar,  Green currents for meromorphic maps of compact K\"ahler manifolds.  {\sl J. of Geometric Analysis},  
arXiv:1107.3063v2 



\item{[BK1]} E. Bedford and KH Kim, On the degree growth of birational mappings in higher dimension,  J. Geom.\  Anal.\ 14 (2004), 567-596. 

\item{[BK2]} E. Bedford and KH Kim,  Periodicities in linear fractional recurrences: Degree growth of birational surface maps, Michigan Math.\  J.\  54 (2006), 647-670.


\item{[BK4]}  E. Bedford and KH Kim,  Linear fractional recurrences: periodicities and integrability, arXiv:0910.4409





\item{[BC]}  J. Blanc and S. Cantat,  personal communication.



\item{[CaL]}      E. Camouzis and G. Ladas,  {\sl Dynamics of Third Order Rational Difference Equations with Open Problems and Conjectures}, Chapman and Hall/CRC Press, 2008.

\item{[C]} S. Cantat,  Dynamique des automorphisms des surfaces $K3$, Acta Math.\ 187 (1): 1--57, 2001. 

\item{[CGM1]}  A. Cima, A. Gasull and F. Ma\~nosas,  On periodic rational difference equations of order $k$. J. Difference Equ. Appl. 10 (2004), no. 6, 549--559.

\item{[CGM2]} A. Cima, A. Gasull and V. Ma\~nosa, Dynamics of the third order Lyness' difference equation,  J. Difference Equ. Appl. 13 (2007), no. 10, 855--884.  

arXiv:math.DS/0612407v1

\item{[CL]} M. Cs\"ornyei and M. Laczkovich, Some periodic and non-periodic recursions, Monatshefte f\"ur Mathematik 132 (2001), 215-236. 

\item{[dFE]}   T. de Fernex   and L. Ein, Resolution of indeterminacy of pairs, {\sl Algebraic geometry}, 165--177, de Gruyter, Berlin, 2002.


\item{[DiF]}  J. Diller and C. Favre, Dynamics of bimeromorphic maps of surfaces, Amer.\ J. of Math., 123 (2001), 1135--1169.
%

\item{[DS]} T.-C. Dinh and N. Sibony, Une borne sup\'erieure pour
l'entropie topologique d'une application rationnelle.   Ann. of Math. (2) 161 (2005), no. 3, 1637--1644.

\item{[DO]}  I. Dolgachev and D. Ortland,  {\sl Point Sets in Projective Spaces and Theta Functions}, Ast\'erisque, Vol.\ 165, 1988.

%

\item{[DuF]}  R. Dujardin and C. Favre,  Distribution of rational maps with a preperiodic critical point. Amer.\  J. Math.\ 130 (2008), no.\ 4, 979--1032.

\item{[EJ]}  A-S Elsenhans and J. Jahnel,  Determinantal quartics and the computation of the Picard group,  arXiv:1010.1923v1

%

%



\item{[GL]}  E.A. Grove and G. Ladas,  {\sl Periodicities  in Nonlinear Difference Equations},  Kluwer Academic Publishers, 2005.

\item{[G1]} V. Guedj, Entropie topologique des applications m\'eromorphes. Ergodic Theory Dynam. Systems 25 (2005), no. 6, 1847Ð1855. 

\item{[G2]} V. Guedj,  Th\'eorie ergodique des transformations rationnelles.   arXiv:math/0611302


\item{[H]} H.P. Hudson,  {\sl Cremona Transformations in Plane and Space}, Cambridge U. Press, 1927.

\item{[KoL]} V.I. Kocic and G. Ladas, {\sl Global Behaviour of Nonlinear Difference Equations of Higher Order with Applications}, Kluwer Academic Publishers 1993.



\item{[KuL]}  M. Kulenovi\'c and G. Ladas,  {\sl Dynamics of Second Order Rational Difference Equations}, CRC Press, 2002.


%

\item{[L]} R.C. Lyness, Notes 1581,1847, and 2952, Math.\ Gazette {\bf 26} (1942), 62, {\bf 29} (1945), 231, and {\bf 45} (1961), 201.

\item{[M1]} C.T. McMullen,   Dynamics on blowups of the projective plane. Publ.\ Math.\ I.H.E.S. No.\ 105 (2007), 49--89.

\item{[M2]}  C.T. McMullen,  K3 surfaces, entropy and glue.  

http://www.math.harvard.edu/$\sim$ctm/papers/home/text/papers/glue/glue.pdf

%



\bigskip

\bigskip
\rightline{Indiana University}

\rightline{Bloomington, IN 47405}

\rightline{\tt bedford@indiana.edu}

\bigskip
\rightline{Florida State University}

\rightline{Tallahassee, FL 32306}

\rightline{\tt kim@math.fsu.edu}

\bye